\documentclass{article}
\pdfoutput=1 

\usepackage{arxiv}

\DeclareMathAlphabet{\pazocal}{OMS}{zplm}{m}{n}
\usepackage{mathrsfs}
\usepackage{calligra}
\usepackage[utf8]{inputenc}
\usepackage[T1]{fontenc} 	
\usepackage{hyperref}
\usepackage{url}            
\usepackage{booktabs}       
\usepackage{amsfonts}       
\usepackage{amssymb,amsmath}
\usepackage{amsthm,titlesec}
\usepackage{graphicx}
\usepackage{enumerate}
\usepackage{enumitem}
\usepackage{bm}
\usepackage{dsfont}
\usepackage{pifont}
\usepackage[super]{nth}

\usepackage{xcolor}
\definecolor{darkpastelgreen}{rgb}{0.01,0.75,0.24}

\usepackage{tabularx}
\usepackage{mathtools}

\theoremstyle{plain}
\newtheorem{theo}{Theorem}[section] 
\newtheorem{lemma}[theo]{Lemma} 
\newtheorem{prop}[theo]{Proposition}
\newtheorem{coroll}[theo]{Corollary}

\theoremstyle{definition}
 
\newtheorem{exmp}[theo]{Example} 
\newtheorem{com}[theo]{Remark}
\newtheorem{dig}[theo]{Digression}

\titleformat{\section}[hang]
  {\bfseries}
  {\thesection.}
  {1ex}
  {}
\titlespacing{\section}{1.5pt}{0.2cm}{0.2cm}

\titleformat{\subsection}[hang]
  {\bfseries}
  {\thesubsection.}
  {1ex}
  {}
\titlespacing{\subsection}{1.5pt}{0.2cm}{0.2cm}

\newcommand{\1}{\partial}
\newcommand{\w}{\omega}
\newcommand{\wn}{\omega_{n}}

\newcommand{\Jxi}{\mathbf{\textit{J}}_{\xi}}

\newcommand{\jota}{\mathscr{J}}

\newcommand{\grad}{\text{grad}_{\textsl{g}}\hspace{0.05cm}}
\newcommand{\dete}{\text{det}\hspace{0.05cm}}
\newcommand{\lna}{\text{log}\hspace{0.05cm} }
\newcommand{\tra}{\text{tr}\hspace{0.05cm} }
\newcommand{\txib}{t_{\textsl{f}}(\xi)}
\newcommand{\txibx}{t_{\textsl{f}}(\xi_x)}
\newcommand{\distP}{\text{dist}_{\textit{p}\hspace{-0.02cm}_{\textit{0}}}}
\newcommand{\pzero}{\textit{p}\hspace{-0.02cm}_{\textit{0}}}
\newcommand{\expo}{\text{exp}}
\newcommand{\hit}{\textit{h}}

\newcommand{\sht}{\textit{\underline{s}}}

\newcommand{\has}{\textsl{H}}
\newcommand{\gamze}{\Gamma_{\hspace{-0.06cm}\textbf{\tiny{0}}}}
\newcommand{\Ricc}{\textsl{R}\hspace{-0.025cm}\text{i\hspace{-0.03cm}c}}

\renewcommand{\shorttitle}{\textit{A SHARP GEOMETRIC INEQUALITY}}

\hypersetup{
pdftitle={A sharp geometric inequality},
pdfsubject={Comparison Geometry, Differential Geometry},
pdfauthor={Adam Rudnik},
pdfkeywords={Asympotically nonnegative curvature, Asymptotic volume ratio, Willmore-type inequality},
}

\begin{document}
\makeatletter
\newcommand{\address}[1]{\gdef\@address{#1}}
\newcommand{\email}[1]{\gdef\@email{\url{#1}}}
\newcommand{\@endstuff}{\par\vspace{\baselineskip}\noindent\small
\begin{tabular}{@{}l}\scshape\@address\\\textit{E-mail address:} \@email\end{tabular}}
\AtEndDocument{\@endstuff}
\makeatother
\title{\boldmath A SHARP GEOMETRIC INEQUALITY FOR CLOSED HYPERSURFACES IN MANIFOLDS WITH ASYMPTOTICALLY NONNEGATIVE CURVATURE}

\author{\href{https://orcid.org/0000-0003-3181-1466}{\includegraphics[scale=0.06]{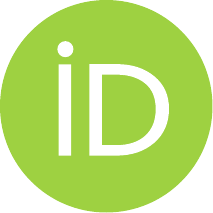}\hspace{1mm}ADAM RUDNIK} 
		}

\address{\parbox{\linewidth}{Department of Mathematics, Institute of Exact Sciences (ICEx), Federal University of Minas Gerais,\\ Belo Horizonte, MG 31270901, Brazil}}
\email{adamrudk@ufmg.br}

\date{}
\maketitle

\begin{abstract}\footnotesize
In this work we establish a sharp geometric inequality for closed hypersurfaces in complete noncompact Riemannian manifolds with asymptotically nonnegative curvature using standard comparison methods in Riemannian Geometry. These methods have been applied in a recent work by Wang \cite{AFST_2023_6_32_1_173_0} to greatly simplify the proof of the Willmore-type inequality in complete noncompact Riemannian manifolds of nonnegative Ricci curvature, which was first proved by Agostiniani, Fagagnolo and Mazzieri \cite{agostiniani2020sharp}.
\end{abstract}


\section{Introduction}

A classical theorem due to Willmore \cite{willmore1968mean} asserts that any compact domain $D \subset \mathbb{R}^{3}$ with regular boundary $\1 D= \Sigma$ has 

\begin{flalign*} 
\hspace{0.8cm} \int_{\Sigma}  \bigg( \frac{\has}{2}\bigg)^2 \textit{d}\sigma & \geq 4\pi&
\end{flalign*}
where $\has$ is the mean curvature of $\Sigma$. Moreover, equality holds only for round spheres. Thus, inside the class of genus zero the inequality above is optimal. It is now widely known that the Willmore inequality extends to $\mathbb{R}^{n}$, $n\geq 3$, to closed submanifolds of any codimension (see e.g. \cite{chen2019geometry}). More generally, in 2020 Agostiniani, Fagagnolo and Mazzieri \cite{agostiniani2020sharp} established a sharp Willmore-type inequality for compact domains in Riemannian manifolds with nonnegative Ricci curvature. Shortly, their inequality states that if $D$ is a compact domain with smooth boundary $\1 D = \Sigma$ in a complete noncompact Riemannian manifold $(\textsl{M},\textsl{g})$ of nonnegative Ricci curvature and $n \geq 3$ then

\begin{flalign*}
\hspace{0.8cm}
\int_{\Sigma} \bigg\vert \frac{\has}{n-1}  \bigg\vert^{n-1} \textit{d}\sigma & \geq \vert \mathbb{S}^{n-1} \vert AVR(\textsl{g}) &
\end{flalign*}
where $AVR(\textsl{g})$ is the asymptotic volume ratio of $(\textsl{M},\textsl{g})$ and $\has$ is the mean curvature of $\Sigma$. They also characterize the equality case provided $AVR(\textsl{g})>0$. The argument in their work is based on PDE methods, and was acknowledged in a subsequent work by Wang \cite{AFST_2023_6_32_1_173_0} as beautiful; nonetheless, highly nontrivial. In his work, Wang obtained the same Willmore-type inequality for noncompact and complete Riemannian manifolds of nonnegative Ricci curvature and his argument is based on standard comparison methods in Riemannian geometry derived from the Riccati equation. \vspace{0.2cm}

In a different direction, Brendle \cite{brendle2021isoperimetric} studied the Isoperimetric problem for minimal hypersurfaces in $\mathbb{R}^{n+1}$ and obtained a Sobolev inequality which holds for any submanifold in Euclidean space. This established a long standing conjecture that states that the isoperimetric inequality holds for compact minimal hypersurfaces in place of domains. We note that Brendle's method, the Alexandrov-Bakelman-Pucci maximum principle, can be extended to more general ambient manifolds as in \cite{brendle2023sobolev} and \cite{johne2021sobolev}. In particular, using Brendle's Sobolev inequalities in nonnegatively curved spaces together with Hölder's inequality is possible to obtain a weaker Willmore-type inequality. Briefly, if $D$ is a compact domain with smooth boundary $\1 D = \Sigma$ in a complete noncompact Riemannian manifold $(\textsl{M},\textsl{g})$ of nonnegative sectional curvature then, as shown in Annex \ref{weaker_willmore_type}, we have

\begin{flalign*}
\hspace{0.8cm}
\int_{\Sigma} \bigg\vert \frac{ \vec{H}}{n-1}  \bigg\vert^{n-1} \textit{d} \sigma & \geq \frac{\vert \mathbb{S}^{n-2} \vert}{n-1} AVR(\textsl{g}) &
\end{flalign*}
while $\alpha_{n-1} > \frac{\alpha_{n-2}}{n-1}$, where $\alpha_{k}$ is the k-dimensional area of the Euclidean unit sphere $\mathbb{S}^{k}$, $\vec{H}$ denotes the mean curvature vector of $\Sigma$, defined in  (\ref{def_mean_curv_vec}) ahead, and $AVR(\textsl{g})$ the asymptotic volume ratio of $(\textsl{M},\textsl{g})$ .\vspace{0.2cm} 

Now we can compare the Willmore's type inequality with the Sobolev inequalities in the above context. On the one hand, the overall extensiveness of these Sobolev inequalities regarding codimension and the flexibility on the function $f$ are broader, not to mention that the claims on these inequalities concern submanifolds possibly with boundary. On the other hand, aside from an application of Hölder's inequality, the constant $AVR(\textsl{g}) \vert \mathbb{S}^{n-1} \vert$ appearing in the Willmore-type inequality is finer, while the hypothesis on the curvature is weaker and the rigidity statement is somehow stronger. \vspace{0.1cm}

Back to the pursuit of extending Brendle's isoperimetric inequality, in 2022, Dong, Lin and Lu \cite{dong2022sobolev} extended Brendle's result to complete noncompact Riemannian manifolds of asymptotically nonnegative curvature. It is important to observe that the theorems obtained by Dong, Li and Lu do not yield our inequality, let aside the curvature assumptions. In particular, the equality cases analyzed in their work force the whole ambient manifold to be isometric with the Euclidean space, while our equality case either bends the ambient to be isometric with a Riemannian manifold of only one end and which is a truncated cone outside an open set $\Omega$, or forces the ambient manifold to be isometric with a specific type of warped product outside an open set, see Theorem \ref{main_theo} for details. Moreover, in the first case, it is exactly the geometry of that set $\Omega$ that possibly prevents the ambient to be isometric with the Euclidean space, see Annex \ref{geometry_of_Omega}. \vspace{0.2cm}

In this work we prove a Willmore-type inequality for the boundary of a compact domain in noncompact and complete Riemannian manifolds with asymptotically nonnegative Ricci curvature, extending the main theorem in \cite{agostiniani2020sharp} to this setting. The class of asymptotically nonnegatively curved Riemannian manifolds was first studied by Abresch in the nineteen eighties, \cite{abresch1985lower} and \cite{abresch1987lower}. In the nineties, Zhu extended the classical Bishop-Gromov inequalities \cite{gromov1999metric} in the case of complete Riemannian manifolds of nonnegative Ricci curvature to these manifolds \cite{10.2307/2374996}. It is interesting to notice that, according to \cite{abresch1985lower}, the class of asymptotically nonnegatively curved spaces includes the class of asymptotically flat manifolds. In particular, it is expected that a Willmore-type inequality for these spaces gives a lower bound for the area of minimal submanifolds, provided that the asymptotic volume ratio does not vanish. \vspace{0.2cm}

Recall that $(\textsl{M},\textsl{g})$ is said to have \textbf{\textit{asymptotically nonnegative}} Ricci (resp. sectional) curvature with base point $\pzero$ and associated function $\lambda$ if $\lambda:[0,\infty) \to \mathbb{R}^{+}$ is a continuous, nonnegative and nonincreasing function such that

\begin{flalign*}
\hspace{0.8cm}
\Ricc_w \geq & -(n-1) (\lambda \circ \distP)(w)\textsl{g}_w
 \hspace{0,5cm}
\big(\text{resp.} \hspace{0.2cm}\text{sec}_{w} \geq - (\lambda \circ \distP)(w)\big),  &
\end{flalign*}
for $w \in  \textsl{M}$ where $\distP(w)$ is the distance from $w$ to $\pzero$ and the associated function $\lambda$ is assumed to satisfy

\begin{flalign}
\label{definition_b_0}
\hspace{0.8cm}
b_0  \doteq & \int_{0}^{\infty} t \lambda(t) dt < \infty &
\end{flalign}
which implies
\begin{flalign}
\label{definition_b_1}
\hspace{0.8cm}
b_1  \doteq & \int_{0}^{\infty}  \lambda(t) dt < \infty. &
\end{flalign}
Note that from the definition of the associated function $\lambda$, either $\lambda$ vanishes identically and we write $\lambda =0$ or there exists $\delta >0$ such that $\lambda >0$ on $[0, \delta )$, in which case $\lambda \neq 0$. In the following, we write $\lambda \neq 0$ meaning that $\lambda$ does not vanish identically. Now, according to Lemma 2.1 in \cite{10.2307/2374996}, the function 

\begin{flalign*}
\hspace{0.8cm}\Theta (r) &= \frac{\textit{vol}(B_{\pzero}(r))}{\w _{n} r^{n}}&
\end{flalign*}
is nonincreasing on $[0, \infty )$ and $0 \leq  \Theta (r) \leq e^{(n-1)b_0}$, where $B_{\pzero}(r)$ is the metric ball of radius $r$ centered at $\pzero$ and $\w_n$ denotes the volume of the $n$-dimensional unit ball in Euclidean space, $\w_n = \vert \mathbb{B}^{n} \vert$. Therefore, we may introduce the asymptotic volume ratio of $(\textsl{M},\textsl{g})$ by

\begin{flalign}
\label{AVRg}
\hspace{0.8cm}
AVR(\textsl{g})& = \lim_{r \rightarrow \infty} \Theta (r) & 
\end{flalign}

Throughout this work we consistently assume that the dimension of the underlying ambient manifold is at least $3$ and that it has positive asymptotic volume ratio, for otherwise our inequalities are trivial.  \vspace{0.2cm}

The next is the main result of this paper. \vspace{0.2cm}

\begin{theo} $\big[$\textit{Willmore-type inequality in asymptotic nonnegatively curved spaces} $\big]$ \newline
\label{main_theo}
Let $(\textsl{M},\textsl{g})$ be a noncompact, complete n-dimensional Riemannian manifold with asymptotically nonnegative Ricci curvature with base point $\pzero$ and associated function $\lambda$. Let $\Omega$ be an open and bounded set with smooth boundary $\Sigma = \partial \Omega$, whose mean curvature is $\has$. Then

\begin{flalign}
\label{willmorex_type_inequality_main}
\hspace{0.8cm}
e^{(n-1)b_0}\int_{\Sigma}\bigg(\bigg \vert \frac{ \has(x)}{n-1} \bigg \vert (1+b_0) + b_1 \bigg)^{n-1}d\sigma (x) & \geq AVR(\textsl{g}) \vert\mathbb{S}^{n-1} \vert &
\end{flalign}
where $AVR(\textsl{g})$ is the asymptotic volume ratio of $\textsl{g}$. In the following, we distinguish two cases.

\begin{enumerate}[label=($\mathbf{\textsl{W}}$\hspace{-0.02cm}\textbf{\arabic*}), leftmargin=1.7cm]
\item If $\lambda = 0$ so that $(\textsl{M},\textsl{g})$ has nonnegative Ricci curvature then equality holds iff $\textsl{M}\setminus \Omega$ is isometric to

\begin{flalign*}
\hspace{0.8cm} \bigg( [r_0, \infty) \times \Sigma, dr\otimes dr + \Big(\frac{r}{r_0}\Big)^2 \textsl{g}_{\Sigma}\bigg), \hspace{0.2cm} \text{where} \hspace{0.2cm} r_{0}^{n-1} &= \bigg(\frac{ \vert \Sigma \vert}{AVR(\textsl{g}) \vert \mathbb{S}^{n-1}\vert}\bigg).&
\end{flalign*}
In particular, $\Sigma$ is a connected totally umbilic hypersurface with constant mean curvature;

\item If $\lambda \neq 0$ and equality holds then $\Sigma$ is a totally umbilic hypersurface with nonnegative constant mean curvature on its connected components. If $\Sigma$ is connected then equality in (\ref{willmorex_type_inequality_main}) holds if and only if $M \setminus \Omega$ is isometric to 

\begin{flalign}
\label{iso_asy_to}
\hspace{0.8cm} \bigg( [r_0, \infty) \times \Sigma & , dr\otimes dr + \varrho(r)^2 \textsl{g}_{\Sigma}\bigg),  \hspace{0.2cm} \text{where} \hspace{0.2cm} r_{0} = \text{dist}_{\textsl{g}}(\pzero, \Sigma),&
\end{flalign}
and $\varrho$ is the solution to the Jacobi equation $\varrho^{\prime \prime}(r) - \lambda(r)\varrho(r)=0$, over the interval $[r_0, \infty)$, with $\varrho(r_0)=1$, $\varrho^{\prime}(r_0)=\frac{\has}{n-1}$ and satisfying the following property:

\begin{flalign}
\label{iso_asy_to}
\hspace{0.8cm} & \lim_{r \rightarrow \infty}\frac{\varrho(r+r_0)}{e^{b_0}\Big(\frac{H}{n-1}(1+b_0) +b_1\Big)r} = 1&
\end{flalign}
In particular, $\Omega$ is a geodesic ball centered at $\pzero$.
\end{enumerate}
\end{theo}

Note that if $b_0 = 0$ in the inequality (\ref{willmorex_type_inequality_main}) above so that $(\textsl{M},\textsl{g})$ has nonnegative Ricci curvature, we obtain 
\begin{flalign}
\label{willmore_type_inequality_main_ricci_nonnegative}
\hspace{0.8cm}
\int_{\Sigma}\bigg \vert \frac{ \has}{n-1} \bigg \vert ^{n-1}d\sigma & \geq AVR(\textsl{g}) \vert\mathbb{S}^{n-1} \vert ,&
\end{flalign}
which is the main result in \cite{agostiniani2020sharp}. Therefore, the inequality (\ref{willmorex_type_inequality_main}) generalizes the Willmore-type inequality for noncompact complete Riemannian manifolds of nonnegative Ricci curvature to asymptotic nonnegatively curved spaces.\vspace{0.2cm}

The issue of the decay of the associated function $\lambda$ is of delicate nature and it is related with various asymptotic behavior of quantities arising from the estimates obtained in our elementary inequalities, e.g. the volume element. In Section \ref{sharp_dom_asym_non} we completely describe this matter and work out a sufficient condition for the inequality (\ref{willmorex_type_inequality_main}) to be strict depending solely on the decay of $\lambda$, see Proposition \ref{prp_decay} for details. \vspace{0.2cm}

The disadvantage of making comparison geometry in manifolds with asymptotically nonnegative curvature is caused by the dependence of the base point $\pzero$. However, because of the decay condition on the curvature, it is expected that the behavior of these manifolds at infinity resemble the behavior of those manifolds with nonnegative curvature. \vspace{0.2cm}

Next, we introduce a class of manifolds that include those with asymptotically nonnegative curvature. This new class will be useful to explore some specific types of warped product spaces to which we can extend the Willmore-type inequality from Theorem \ref{main_theo}. Moreover, we investigate to what extent we may bring off the rigidity statement in the main theorem, not only for null-homologous hypersurfaces but also for a broader class of hypersurfaces. \vspace{0.2cm}

Let $(\textsl{M},\textsl{g})$ be a complete noncompact $n$-dimensional Riemannian manifold. We say that $\textsl{g}$ has \textbf{\textit{asymptotically$^{\text{\ding{93}}}$ nonnegative}} Ricci (resp. sectional) curvature relatively to a $k$-dimensional submanifold $\textsl{N}^{k} \subset \textsl{M}$ with associated function $\lambda$ if $\lambda:[0,\infty) \to [0,\infty)$ is a continuous, nonnegative and nonincreasing function satisfying (\ref{definition_b_0}) and such that

\begin{flalign*}
\hspace{0.8cm}
\Ricc_w \geq & -(n-1) (\lambda \circ \text{dist}_{\textsl{N}})(w)\textsl{g}_w
 \hspace{0,5cm}
\big(\text{resp.} \hspace{0.2cm} \text{sec}_{w} \geq - (\lambda \circ \text{dist}_{\textsl{N}})(w)\big), &
\end{flalign*}
for $w \in  \textsl{M}$, where $\text{dist}_{\textsl{N}}(w)$ is the distance from $w$ to $\textsl{N}$. Observe that if $\textsl{N}$ is a single point $\pzero$, then the idea of asymptotic$^{\text{\ding{93}}}$ nonnegative curvature reduces to the concept of asymptotic nonnegative curvature. In addition, it is not clear whether asymptotic$^{\text{\ding{93}}}$ nonnegatively curved spaces possess a well defined notion of asymptotic volume ratio. We address this issue at Section \ref{new} and restrict ourselves to spaces that have this property as defined there.\vspace{0.2cm}

\begin{theo}  $\big[$\textit{Willmore-type inequality in asymptotic$^{\text{\ding{93}}}$ nonnegatively curved spaces} $\big]$ \newline
\label{sec_main_theo}
Let $(\textsl{M},\textsl{g})$ be a noncompact, complete	$n$-dimensional Riemannian manifold with asymptotically$^{\text{\ding{93}}}$ nonnegative Ricci curvature relatively to a closed and connected hypersurface $\gamze$ and associated function $\lambda$, with $\lambda \neq 0$. Suppose that $\gamze$ separates $\textsl{M}$ into two connected components: $\textsl{M} \setminus \gamze = \textsl{M}^{+} \cup \textsl{M}^{-}$, with $\textsl{M}^{+}$ unbounded. Let $\Sigma \subset \textsl{M}^{+}$ be a closed hypersurface homologous to $\gamze$. Let $\Omega \subset \textsl{M}$ be the open set with the property $\1 \Omega = \gamze \cup \Sigma$ and assume further that $\vert \Omega \vert < \infty$. Then

\begin{flalign}
\label{wil_typ_inq_sec_main}
\hspace{0.8cm}
e^{(n-1)b_0}\int_{\Sigma}\bigg(\bigg \vert \frac{ \has(x)}{n-1} \bigg \vert (1+b_0) + b_1 \bigg)^{n-1}d\sigma (x) & \geq AVR^{\text{\ding{93}}}(\textsl{g}) \vert\mathbb{S}^{n-1} \vert &
\end{flalign}
where $\has$ is the mean curvature of $\Sigma$ and $AVR^{\text{\ding{93}}}(\textsl{g})$ is the asymptotic$^{\text{\ding{93}}}$ volume ratio of $\textsl{g}$ as defined in Digression \ref{asp_vol_rti_str}. Moreover, if equality holds then $\Sigma$ is a totally umbilic hypersurface whose mean curvature is constant on its connected components. In addition to that, if $\Sigma$ is connected then equality in (\ref{wil_typ_inq_sec_main}) holds iff $\textsl{M}^{+} \setminus \Omega$ is isometric to 

\begin{flalign}
\label{iso_asy_star_to}
\hspace{0.8cm} \bigg( [r_0, \infty) \times \Sigma & , dr\otimes dr + \varrho(r)^2 \textsl{g}_{\Sigma}\bigg),  \hspace{0.2cm} \text{where} \hspace{0.2cm} r_{0} = \text{dist}_{\textsl{g}}(\gamze, \Sigma),&
\end{flalign}
and $\varrho$ is the solution to the Jacobi equation $\varrho^{\prime \prime}(r) - \lambda(r)\varrho(r)=0$, over $[r_0, \infty)$, with $\varrho(r_0)=1$, $\varrho^{\prime}(r_0)=\frac{\has}{n-1}$ and obeying to

\begin{flalign}
\label{iso_asy_to_star}
\hspace{0.8cm} & \lim_{r \rightarrow \infty}\frac{\varrho(r)}{e^{b_0}\Big(\frac{H}{n-1}(1+b_0) +b_1\Big)r} = 1.&
\end{flalign}
In particular, $\gamze$ is equidistant to $\Sigma$.
\end{theo}

As an immediate application of our inequality we have the following estimate. \vspace{0.3cm}

\begin{coroll}  $\big[$\textit{Lower bound for the area of closed minimal hypersurfaces} $\big]$ \newline
\label{area_lower_bound_minimal_hypersurfaces}
As in Theorem \ref{main_theo} (resp. Theorem \ref{sec_main_theo}) let $\Sigma$ be a closed hypersurface in a complete noncompact n-dimensional Riemannian ambient $(\textsl{M},\textsl{g})$ with asymptotically (resp. asymptotically$^{\text{\ding{93}}}$) nonnegative Ricci curvature with base point $\pzero$ (resp. relatively to $\gamze$) and associated function $\lambda$, with $\lambda \neq 0$. If $\Sigma$ is minimal and null-homologous (resp. homologous to $\gamze$) then,

\begin{flalign}
\label{lower_bound_minimal_inequality}
\hspace{0.8cm}
\vert \Sigma \vert & \geq \frac{\mathcal{A}(\textsl{g})\vert\mathbb{S}^{n-1} \vert}{\hspace{0.2cm}(e^{b_0}b_{1})^{n-1}}  ,&
\end{flalign}
where $\mathcal{A}(\textsl{g})$ is the asymptotic (resp. asymptotic$^{\text{\ding{93}}}$) volume ratio of $\textsl{g}$. \vspace{0.2cm}
\end{coroll}

There is a nice geometrical idea behind the lower bound above. If $b_1$ approaches zero then the RHS in (\ref{lower_bound_minimal_inequality}) becomes very large. Also, the associated function $\lambda$ is approaching zero and in doing so, the whole Ricci tensor is gradually changing from being asymptotically nonnegative to being nonnegative. However, noncompact spaces of nonnegative Ricci curvature with Euclidean volume growth do not support closed minimal hypersurfaces, so it is expected that the area of $\Sigma$ is getting larger, as formula (\ref{lower_bound_minimal_inequality}) tells. \vspace{0.2cm}

It is easy to check that if $(\textsl{M},\textsl{g})$ has asymptotically nonnegative Ricci curvature with base point $\pzero$ then $(\textsl{M},\textsl{g})$ has asymptotically$^{\text{\ding{93}}}$ nonnegative Ricci curvature relatively to $\textsl{K}$, for any compact submanifold $\textsl{K}$ such that $\pzero \in \textsl{K}$. In this line of thought, the question whether there are more general volume comparison theorems of Heintze-Karcher type in asymptotically nonnegative curved spaces remains open. The case where $\textsl{N} =\{ \pzero\}$ and the case where $\textsl{N}$ is a closed geodesic were already inspected (see the remark after Corollary 2.1 in \cite{10.2307/2374996}). In the present work we treat the codimension one, closed and embedded case. \vspace{0.2cm}

Next we apply our inequality (\ref{wil_typ_inq_sec_main}) for a class of warped product manifolds that we now describe. Let $(\textsl{N}, \textsl{g}_{\textsl{N}})$ be a connected and compact Riemannian manifold of dimension $n-1 \geq 2$ such that 

\begin{flalign*}
\hspace{0.8cm}
\Ricc_{\textsl{N}} & \geq (n-2)\rho  \textsl{g}_{\textsl{N}}, &
\end{flalign*}
for some constant $\rho$. Furthermore, consider a positive and smooth function $h:[0,\infty) \to \mathbb{R}$ and define the quantities $\lambda_1, \lambda_2:[0, \infty) \to \mathbb{R}$ by

\begin{flalign*}
\hspace{0.8cm}
\lambda_1 (r) & = \frac{h^{\prime \prime}(r)}{h(r)}\hspace{0.2cm} \text{and}\hspace{0.2cm}  \lambda_2(r) = \frac{1}{(n-1)}\frac{h^{\prime \prime}(r)}{h(r)} -\bigg(\frac{n-2}{n-1}\bigg)\frac{\rho - h^{\prime}(r)^2}{h(r)^2}. &
\end{flalign*}
In terms of the functions $\lambda_1$ and $\lambda_2$, assume that the function $h$ satisfies the following conditions:

\begin{enumerate}[label=($\mathbf{\Lambda}$\hspace{-0.02cm}\textbf{\arabic*}),leftmargin=1.7cm]
\item Either we have $\lambda_1 > 0$ or $\lambda_2 > 0$ at some point.

\item Assume that there exists a continuous monotone decreasing nonnegative function $\lambda:[0, \infty)\to \mathbb{R}$ such that, for all $r \in [0,\infty)$, we have $\lambda(r) \geq \text{max}\{ \lambda_1(r), \lambda_2(r) \}$ and 

\begin{flalign}
\hspace{0.8cm}
b_0 \doteq & \int_{0}^{\infty}t \lambda(t)dt < \infty. &
\end{flalign}
\item The function $[0, \infty) \ni r \mapsto \frac{h(r)}{r}$ is eventually nonincreasing.
\item There exists $\tau_0 \geq 0$ such that $h^{\prime}(r) \geq 0$ for all $r \in [\tau_0, \infty)$ and $\displaystyle{\lim_{r \rightarrow \infty}h(r)= \infty}$. \vspace{0.2cm}
\end{enumerate}

We now consider the product $\textsl{M} \doteq [0, \infty)\times \textsl{N}$ endowed with the Riemannian metric 

\begin{flalign}
\label{warped_product_metric}
\hspace{0.8cm}
\textsl{g} & \doteq dr \otimes dr + h(r)^2 \textsl{g}_{N}. &
\end{flalign}

Our main application of Theorem \ref{sec_main_theo} is the following. \vspace{0.2cm}

\begin{coroll} $\big[$\textit{Application to warped product spaces} $\big]$ \newline
\label{main_prop}
Let $(\textsl{M},\textsl{g})$ be a warped product manifold satisfying conditions $(\mathbf{\Lambda}\hspace{-0.02cm}\textbf{1})\textit{---}(\mathbf{\Lambda}\hspace{-0.02cm}\textbf{3})$. Then, $\gamze = \{0\} \times \textsl{N}$ is a compact hypersurface whose area satisfies

\begin{flalign}
\label{willmorex_type_inequality_wp}
\hspace{0.8cm}
\big \vert \gamze \big \vert & \geq \frac{AVR^{\text{\ding{93}}}(\textsl{g}) \vert\mathbb{S}^{n-1} \vert}{\hspace{0.15cm}\bigg( e^{b_0}\Big( \Big\vert\frac{h^{\prime}(0)}{h(0)}\Big\vert (1+b_0 )+ b_1\Big)\bigg)^{n-1}}. &
\end{flalign}
Let $\Sigma \subset \textsl{M}$ be a connected and closed hypersurface with mean curvature $\has$. If $\Sigma$ is homologous to $\gamze$ then

\begin{flalign}
\label{will_type_ine_connected}
\hspace{0.8cm}
e^{(n-1)b_0}\int_{\Sigma}\bigg(\bigg \vert \frac{ \has(x)}{n-1} \bigg \vert (1+b_0) + b_1 \bigg)^{n-1}d\sigma (x) & \geq AVR^{\text{\ding{93}}}(\textsl{g}) \vert\mathbb{S}^{n-1} \vert . &
\end{flalign}
Equality in (\ref{will_type_ine_connected}) holds iff $\Sigma$ is a slice $\{r_0\}\times \textsl{N}$, $h$ is the solution to $h^{\prime \prime}(r)- \lambda(r)h(r)=0$ over $[r_0, \infty)$, $h(r_0)=1$, $h^{\prime}(r_0)=\has/(n-1)$ and $\lim_{r \rightarrow \infty} h(r)/\big[ e^{b_0}\big(h^{\prime}(r_0)(1+b_0) +b_1\big)r\big]=1$.

\end{coroll}
If hypothesis ($\mathbf{\Lambda}$\hspace{-0.02cm}\textbf{1}) does not hold then $(\textsl{M},\textsl{g})$ has $\Ricc_{\textsl{M}} \geq 0$. Then,
hypothesis ($\mathbf{\Lambda}$\hspace{-0.02cm}\textbf{1}) suggests that the ambient does not have, possibly, nonnegative Ricci curvature. For instance, if $(\textsl{N}, \textsl{g}_{\textsl{N}})$ is Einstein, meaning that $\Ricc_{\textsl{N}} = (n-2)\rho  \textsl{g}_{\textsl{N}}$ then $(\textsl{M}, \textsl{g})$ does not have nonnegative Ricci curvature. The hypothesis ($\mathbf{\Lambda}$\hspace{-0.02cm}\textbf{2}) guarantees that $(\textsl{M}, \textsl{g})$ has asymptotically$^{\text{\ding{93}}}$ nonnegative Ricci curvature relatively to $\gamze$. The condition ($\mathbf{\Lambda}$\hspace{-0.02cm}\textbf{3}) is there to ensure finiteness of the asymptotic$^{\text{\ding{93}}}$ volume ratio of $\textsl{g}$ and, finally, ($\mathbf{\Lambda}$\hspace{-0.02cm}\textbf{4}) guarantees that $AVR^{\text{\ding{93}}}(\textsl{g})$ is positive. See Lemma \ref{bas_prp_mg} for details. \vspace{0.2cm}

We point out that the Schwarzschild and the Reissner-Nordstrom manifolds meet the hypothesis of Corollary \ref{main_prop}. We briefly recall the definition of these spaces here. Let $m>0$ be a real number and write $\{ s>0: 1-ms^{2-n} \} = (\sht, \infty)$. The Schwarzschild manifold is defined by 

\begin{flalign}
\label{n-d_sch_mtr}
\hspace{0.8cm}
 \textsl{M}_{\textit{S}} & = (\sht, \infty) \times \mathbb{S}^{n-1}, \hspace{0.2cm} \textsl{g}_{\textit{S}} = \frac{1}{1-ms^{2-n}}ds \otimes ds + s^2 \textsl{g}_{\mathbb{S}^{n-1}}.&
\end{flalign}
In Section \ref{app_und_exm} we discuss in depth this example: we show how to apply the inequality (\ref{willmorex_type_inequality_wp}) to the horizon of the Schwarzschild space. The Reissner-Nordstrom manifold is defined by 

\begin{flalign}
\label{n-d_rei_nord_mtr}
\hspace{0.8cm}
\textsl{M}_{\textit{R}} & = (\sht, \infty) \times \mathbb{S}^{n-1}, \hspace{0.2cm} \textsl{g}_{\textit{R}} = \frac{1}{1-ms^{2-n}+q^2s^{4-2n}}ds \otimes ds + s^2 \textsl{g}_{\mathbb{S}^{n-1}},&
\end{flalign}
where $m > 2q>0$ are constants and $\sht$ is defined as the larger of the $2$ solutions of the equation $1-ms^{2-n}+q^2s^{4-2n} = 0$. \vspace{0.2cm}

We remark that the simplification in \cite{AFST_2023_6_32_1_173_0} was foundational to our inequalities, but the comparison techniques used in our work are based on the Jacobi equation while those used in \cite{AFST_2023_6_32_1_173_0} are based on the Riccati equation which, in general, works better with comparison theory. It is also important to mention that, regarding the elementary inequalities derived from the Jacobi equation, our techniques are based in \cite{10.2307/2374996}. Moreover, concerning the applications to warped product spaces and the examples described here, the authors were influenced by \cite{brendle2013constant}. Finally, our proof shares some common features with the work of Heintze and Karcher \cite{heintze1978general} regarding the volume of tubular neighborhoods. \vspace{0.2cm}

This paper is organized as follows. Section \ref{notions_und_notations} establishes the preliminary tools for the remaining of the text. In Section \ref{section_elel_ineq} we show some elementary inequalities that will play a pivotal role in this work. These inequalities are based on ODE methods that go back to Sturm, but which are quite trivial whence no reference is given. The proof of the main theorem, Theorem \ref{main_theo}, occupies Section \ref{sharp_dom_asym_non} and it is divided into two steps: the main inequality and the discussion of the rigidity statement. In Section \ref{new}, we develop the adaptation of the proof of the main theorem to the class of asymptotic$^{\text{\ding{93}}}$ nonnegatively curved spaces, namely Theorem \ref{sec_main_theo}. Next, Section \ref{app_und_exm} is reserved to discuss how Corollary \ref{main_prop} may be obtained from Theorem \ref{sec_main_theo}, we postpone the converse of Theorems \ref{main_theo} and \ref{sec_main_theo} to this section and give concrete examples of our inequality. Finally, in Section \ref{zeta} we conclude this work discussing in detail some observations mentioned throughout the text.

\section{Notions and notations}
\label{notions_und_notations}

Let $\Sigma$ be, once and for all, a fixed embedded and closed hypersurface in $(\textsl{M},\textsl{g})$. Denote by $w$ the variable in $\textsl{M}$ and by the letter $x$ the variable in $\Sigma$. We either take $\Sigma$ to be null-homologous, in which case $\1 \Omega = \Sigma$, for some open and bounded set $\Omega$, or we take it to be homologous to some connected hypersurface $\gamze$ in which case we write $\1 \Omega = \Sigma \cup \gamze$, for some open set $\Omega$. In the latter, we also assume that $\textsl{M} \setminus \gamze = \textsl{M}^{+} \cup \textsl{M}^{-}$ where $\textsl{M}^{+}$ and $\textsl{M}^{-}$ are disjoint open subsets with $\Sigma \subset \textsl{M}^{+}$. Either way, we consider the oriented distance

\begin{flalign}
\label{distance_to_hyper}
\hspace{0.8cm}
r(w)& = \text{dist}_{\textsl{g}}(w, \Sigma), \hspace{0.2cm} w \in \textsl{M}\setminus \Omega \hspace{0.2cm} (\text{or} \hspace{0.2cm} w \in \textsl{M}^{+}\setminus \Omega). &
\end{flalign}
Let $\xi :\Sigma \to T\textsl{M}$ be the outward-pointing (relatively to $\Omega$) unit, smooth and normal to $\Sigma$ vector field. Let $\nu \Sigma$ be the \textbf{\textit{normal bundle}} of the hypersurface $\Sigma$. The \textbf{\textit{normal exponential map}} is the restriction of the exponential map of $\textsl{M}$ to $\nu \Sigma$ and it is denoted by $\expo^{\perp}$. Indicate by $S\nu \Sigma$ the \textbf{\textit{unit normal bundle}}, that is, those vectors in $\nu \Sigma$ of norm $\mathds{1}$. Given $x \in \Sigma$, consider the geodesic $\gamma_{\xi_x} (t) = \expo_{x}t\xi_x $. Define the \textbf{\textit{cut time of $(x, \xi_x)$}} by 

\begin{flalign}
\label{cut_time}
\hspace{0.8cm} \tau_{c}(\xi_x)& \doteq \text{sup } \{ b>0: \hspace{0.1cm} r(\gamma_{\xi_x}(b)) = b  \} \in (0,\infty]&
\end{flalign}
where $r=r(w)$ is the distance (\ref{distance_to_hyper}). We say that the geodesic $\gamma_{\xi_x}$ realizes the distance in the interval $[0, \tau_{c}(\xi_x)]$. The \textbf{\textit{cut locus of $\Sigma$}}, denoted by $C( \Sigma)$, is the image of $\{ \tau_c(\zeta)\zeta : \zeta \in S\nu \Sigma \hspace{0.2cm} \text{and} \hspace{0.2cm} \tau_c{(\zeta)} < \infty \}$ under the normal exponential map. As in the case where $\Sigma$ is a point, the function $\tau_c: S\nu \Sigma \to \mathbb{R}$ is continuous. Since, for each $x \in \Sigma$, there is a unique unit outward-pointing normal vector $\xi_x$ we can write $\tau_c(x)$ or $\tau_{c}(\xi_x)$ interchangeably and consider $\tau_{c}$ as a continuous function over $\Sigma$. Furthermore, the set $C(\Sigma)$ is closed and the normal exponential map

\begin{flalign*}
\hspace{0.8cm}
\expo^{\perp} : &\{ t\xi_x : x \in \Sigma \hspace{0.2cm} \text{and} \hspace{0.2cm} 0 \leq t < \tau_{c} (x)\} \to \mathcal{U} \setminus C(\Sigma)&
\end{flalign*}
is a diffeomorphism, where $\mathcal{U} = \textsl{M}\setminus \Omega$ or $\mathcal{U} = \textsl{M}^{+} \setminus \Omega$, depending on whether $\Sigma$ is null-homologous or not. These assertions follow from standard arguments (see e.g. \cite{petersen2006riemannian}, pp. 139 \textit{--} 141). As a result, there exists a diffeomorphism 

\begin{flalign}
\label{diffeo_Phi}
\hspace{0.8cm}
\Phi &: E  \doteq \{ (x,r) \in \Sigma \times [0, \infty) : r < \tau_c(x)  \}  \to \mathcal{U} \setminus C( \Sigma), \hspace{0.2cm} \Phi(x, r) = \expo_{x}r\xi_x.  &
\end{flalign}
We will routinely consider the pullback of the volume element in $\textsl{M}$ induced by $\textsl{g}$ to $E$ via this diffeomorphism, that is, $(\Phi^{*}d \textit{vol}_{\textsl{g}})_{(x, r)} = \dete J_{\xi_x}(r) d\sigma \wedge dr$, for $(x, r) \in E$. Opportunely, using Jacobi Fields, we will give a description of the coefficient function $\dete J_{\xi}(t)$. For more details, refer to \cite{ballmann2016riccati}.

The symmetric bilinear \textbf{\textit{second fundamental form of}} $\Sigma$ is the $\binom{2}{1}$-tensor field defined by 

\begin{flalign*}
\hspace{0.8cm}
\alpha(X,Y) & = (\nabla_{X}Y)^{\perp}, \hspace{0.2cm} X,Y \in \mathfrak{X}(\Sigma) &
\end{flalign*}
and the tensor field of self-adjoint endomorphisms $A_{\xi}$ of $\Sigma$ defined by $\langle A_{\xi}X,Y \rangle = \langle \alpha( X,Y), \xi \rangle $, where $X,Y \in \mathfrak{X}(\Sigma)$, is called the \textbf{\textit{shape operator of $\Sigma$}}. The \textbf{\textit{mean curvature vector $\vec{H}$ of $\Sigma$}} at $x$ is the normal vector defined by

\begin{flalign}
\label{def_mean_curv_vec}
\hspace{0.8cm}
\vec{H}(x) & \doteq \sum_{i=1}^{n-1}\alpha (X_i, X_i) &
\end{flalign}
in terms of an orthonormal basis $X_1, ..., X_{n-1}$ of $T_{x}\Sigma$. \vspace{0.1cm}

For $t> 0$, we say that $\gamma_{\xi}(t)$ is a \textbf{\textit{focal point}} of $\Sigma$ along the geodesic $\gamma_{\xi}$ if there is a nontrivial Jacobi field $J$ along $\gamma_{\xi}$ satisfying \textit{(i)} $J(0) \in T_{p}\Sigma $, \textit{(ii)} $D_tJ(0) + A_{\xi}J(0) \in  T_{p}\Sigma^{\perp}$ and \textit{(iii)} $J(t)=0$, where $D_tJ$ is the usual covariant differentiation of $J$ along $\gamma_{\xi}$. If there is such a $t>0$, then we let $\txib$ be the smallest such $t$ and call $\gamma_{\xi}(t_{\textsl{f}}(\xi))$ the first focal point of $\Sigma$ along $\gamma_{\xi}$. We have $\tau_{c}(\xi) \leq t_{\textsl{f}}(\xi)$, as in the case when $\Sigma$ is a point. Notice that we may have $t_{\textsl{f}}(\xi)=\infty$. For more details on the matter refer to \cite{do1992riemannian} (pp. 227 \textit{--} 235). \vspace{0.1cm}

Denote the space of vector fields along $\gamma = \gamma_{\xi}$ that are normal to $\dot{\gamma}$ by $\mathcal{J}^{\perp}(\gamma)$. A \textbf{\textit{Jacobi tensor field}} $\mathbf{J}$ along $\gamma$ is a smooth $\binom{1}{1}$-tensor field of endomorphisms $\mathcal{J}^{\perp}(\dot{\gamma}) \to \mathcal{J}^{\perp}(\dot{\gamma})$ which satisfies the Jacobi equation

\begin{flalign*}
\hspace{0.8cm}
\mathbf{J}^{\prime \prime}& + \textsl{R}_{\dot{\gamma}}(\mathbf{J}) = 0 &
\end{flalign*}
where $\textsl{R}_{\dot{\gamma}}$ is the endomorphism $X \mapsto \textsl{R}_{\dot{\gamma}}(X) = \textsl{R}(X, \dot{\gamma}) \dot{\gamma}$. The following construction will be used throughout the text: let $(e_1,...,e_{n-1}) \subset T_{x}\Sigma$ be an orthonormal basis consisting of eigenvectors of the shape operator $A_{\xi}$ and let $(E_{\alpha})$ be a parallel $(n-1)$ orthonormal frame along the geodesic $\gamma$ such that $E_{\alpha}(0) = e_{\alpha}$. We then consider the operator $\Jxi = \Jxi(t)$ which takes the velocity $\dot{\gamma}$ to itself and maps $\mathcal{J}^{\perp}(\gamma)$ into $\mathcal{J}^{\perp}(\gamma)$ according to the rule

\begin{flalign}
\label{definition_operator_Jac}
\hspace{0.8cm} \Jxi(t)E_{k}(t)& = J_k(t), \hspace{0.2cm}k=1,...,n-1  &
\end{flalign}
where $J_k$ is the Jacobi vector field with initial conditions $J_{k}(0)=  e_k$ and $D_t J_{k}(0) = \lambda_{k}e_{k}$, where the $\lambda_{k}'s$ fulfill $A_{\xi}e_k = \lambda_k e_k$. Then, $\Jxi$ is easily seen to be a Jacobi tensor field along $\gamma$. Moreover, the field of endomorphisms $U(t)  \doteq \Jxi^{\prime} \Jxi^{-1}(t)$, which is well defined over the interval $[0, t_{\textsl{f}}(\xi))$, satisfies $UJ_k = D_tJ_k$, $k=1,..., n-1$ and is a solution to the Riccati equation

\begin{flalign}
\label{the_riccati_equation}
\hspace{0.8cm} D_t U& + U^{2}+ \textsl{R}_{\dot{\gamma}}= 0.&
\end{flalign}
In this work, we will be frequently producing quantities relatively to the distance (\ref{distance_to_hyper}). One of the most important tensors associated to a function is its Hessian. In particular, it is well known that the Hessian operator of any local distance function satisfies the Riccati equation (\ref{the_riccati_equation}). It is easy to see that if $\mathcal{H}_{r}$ is the Hessian operator of (\ref{distance_to_hyper}), i.e., $D^{2}r (\cdot, \cdot)= \langle \mathcal{H}_{r} \cdot, \cdot \rangle$, where $D^{2}r$ is the usual $\binom{2}{0}$ Hessian, then $\mathcal{H}_{r} = U$ along the integral curves of $\grad r$. We define the \textbf{\textit{mean curvature of $\Sigma$}} at $x$ by $\has(x) \doteq \tra U(0)$, or in other words, $\has(x) = \Delta r(0)$. By (\ref{def_mean_curv_vec}), we have

\begin{flalign*}
\hspace{0.8cm}
\vec{H} &= - \tra U(0) \1_r = -\has \1_r &
\end{flalign*}
We say that $\Sigma$ has \textbf{\textit{mean convex boundary}} (MCB) if the mean curvature vector $\vec{H}$ is an inward pointing normal. The following example is a particular case of Example 4.3 in \cite{ballmann2016riccati}. \vspace{0.2cm}

\begin{exmp} 
\label{example1}
If $U = \Jxi^{\prime} \Jxi^{-1}$, which is defined as long as $\dete \Jxi >0$, then set

\begin{flalign}
\hspace{0.8cm} u(t) &\doteq \frac{1}{n-1}\tra U(t). &
\end{flalign}
Since $U$ satisfies (\ref{the_riccati_equation}) and covariant differentiation commutes with contractions, we have

\begin{flalign*}
\hspace{0.8cm} 
	 u^{\prime} 
	 & = \frac{1}{n-1} \tra D_t U \notag \\
	 & = -\frac{1}{n-1} \tra U^{2} - \frac{1}{n-1} \tra \textsl{R}_{\dot{\gamma}} \notag \\
	 & \leq - \frac{1}{(n-1)^{2}}\big( \tra U\big)^{2} - \frac{1}{n-1} \tra \textsl{R}_{\dot{\gamma}} \notag \\
	 & = - u^2 - \frac{1}{n-1} \Ricc(\dot{\gamma}, \dot{\gamma}), &
\end{flalign*}
where we have used Schwarz inequality $(\tra U)^2 \leq (n-1) \tra (U^{2})$. Now, if $\Ricc(\dot{\gamma}, \dot{\gamma}) \geq (n-1)\kappa$, where $\kappa(t)$ is a smooth function defined on the domain of definition of $U$ then 

\begin{flalign}
\label{eq_ric_eq}
\hspace{0.8cm} 
u^{\prime}(t)& \leq   - u(t)^{2} - \kappa (t) .&
\end{flalign}
Equality in (\ref{eq_ric_eq}) holds if and only if $U = u \mathcal{P}_{\hspace{-0.08cm}\dot{\gamma}}$ and  $\Ricc(\dot{\gamma}, \dot{\gamma}) = (n-1)\kappa $ along $\gamma$, where $\mathcal{P}_{\hspace{-0.08cm}\dot{\gamma}}$ is the orthogonal projection onto $(\dot{\gamma})^{\perp}$. Since $U$ satisfies the Riccati equation, equality in (\ref{eq_ric_eq}) holds if and only if $U = u \mathcal{P}_{\hspace{-0.08cm}\dot{\gamma}}$ and $ \textsl{R}_{\dot{\gamma}}= \kappa  \mathcal{P}_{\hspace{-0.08cm}\dot{\gamma}}$. $\hspace{0.3cm} \square$
\end{exmp}

In this article we will frequently work with \textit{volumen} of oriented tubular neighborhoods of hypersurfaces. Oriented meaning that a side is chosen, which we assume, by convention, to be the one given by the vector field $\xi$. That being said, if we denote this oriented neighborhood by $\mathcal{T}_{\Sigma}^{+}(r)$, which is the set $\{w \in \mathcal{U}: \text{dist}_{\textsl{g}}(w, \Sigma ) <r \}$, then

\begin{flalign}
\label{avrg_makes_sense_domain}
\hspace{0.8cm} 
\textit{vol}(\mathcal{T}_{\Sigma}^{+}(r)) & = \int_{\Sigma}\int_{0}^{r \wedge \tau_{c}(x)} \dete J_{\xi_x}(r)dr d\sigma(x),&
\end{flalign}
where $r \wedge \tau_{c}(x)$ is the minimum between $r$ and $\tau_c(x)$. A word about the asymptotic volume ratio is in order. Given a compact domain $D$ in a Riemannian manifold $(\textsl{M},\textsl{g})$ with asymptotically nonnegative curvature (Ricci or sectional) with base point $\pzero$, denote the geodesic tube of radius $r$ around $D$ by $\mathcal{T}_{D}(r)$ so that $\textit{vol}(\mathcal{T}_{D}(r)) = \vert D \vert + \textit{vol}(\mathcal{T}_{\1 D}^{+}(r))$. Clearly, if $\pzero \in D$ we have $\Theta (r) \leq \frac{\textit{vol}(\mathcal{T}_{D}(r))}{ \w _{n}r^n}$, for all $r>0$, so that

\begin{flalign}
\label{avrg_makes_sense_domain}
\hspace{0.8cm} AVR(\textsl{g}) &\leq \lim_{r \rightarrow \infty}
\frac{\textit{vol}(\mathcal{T}_{D}(r))}{\w _{n}r^n}. &
\end{flalign}
It is not difficult to see that if $\pzero \notin D$ then, eventually, for sufficiently big radius $r$ the set $B_{\pzero}(r) \setminus \mathcal{T}_{D}(r)$ has small volume in comparison with $r^{n}$ and, therefore, (\ref{avrg_makes_sense_domain}) also holds in this case. \vspace{0.2cm}

\section{Elementary inequalities}
\label{section_elel_ineq}

Throughout this section we assume that $\lambda$ is nontrivial, i.e., $\lambda \neq 0$ and we fix a point $x \in \Sigma$ and a unit speed geodesic $ \gamma_{\xi}$ that starts at $x$ and realizes the distance (\ref{distance_to_hyper}), that is, $r(\gamma_{\xi}(t))=t$ for $t$ in some interval $[0, \alpha)$, where $\alpha >0$. Hence, all quantities considered here depend only on the variable $t$ and the point $x$ is momentarily ignored. Also, from now on, by abuse of notation, the function $\lambda$ is meant to be $(\lambda \circ \distP \circ \gamma_{\xi} )(t).$ We observe that if $H(x)\neq 0$ then all comparisons proved in this section are easier. \vspace{0.1cm}

Let $\Jxi(t)$ be the Jacobi tensor field defined in (\ref{definition_operator_Jac}) along the geodesic $\gamma_{\xi}$ which, from now on, we denote simply by $\gamma$. Consider the function $\mathscr{J} : [0, t_{\textsl{f}}(\xi)) \to \mathbb{R}$ defined by
\begin{flalign}
\label{definition_of_J}
\hspace{0.8cm}
\mathscr{J}(t) & \doteq \big(\dete \Jxi(t)\big)^{1/(n-1)}&
\end{flalign}
and compare it with 4.B.4 in \cite{gallot1990riemannian}. Note that we may always arrange the basis $(e_1,...,e_{n-1})$ of $T_x \Sigma$ so that $\Jxi(0) =$ the identity of $T_x \Sigma$. \vspace{0.2cm}

\begin{lemma}
\label{lemma1_elel} Assume $ \Ricc(\dot{\gamma}(t), \dot{\gamma}(t)) \geq -(n-1) (\lambda \circ \distP)(\gamma (t))$ for all $t$ in the domain of definition of $\gamma$ and note that $\dete \Jxi >0$ on $[0,t_{\textsl{f}}(\xi))$. Then, the function $\jota$ defined in (\ref{definition_of_J}) satisfies
\begin{flalign}
\label{inequalities_satisfies_J}
	\hspace{0.8cm} &  
	\begin{cases}
 	\jota^{\prime \prime} -\lambda \jota \leq 0; & \\
 	\jota(0)=1, \hspace{0.3cm} \jota^{\prime}(0)=\has/(n-1). \hspace{3cm} & 
	\end{cases}&
\end{flalign}
where $\has = \has(x)$ is the mean curvature of $\Sigma$ at $x$. In addition to that, if $U \doteq \Jxi^{\prime} \Jxi^{-1}$ then the equality $\jota^{\prime \prime} -\lambda \jota = 0$ holds if and only if $U = \frac{\tra U}{n-1} \mathcal{P}_{\hspace{-0.08cm}\dot{\gamma}}$ and  $\textsl{R}_{\dot{\gamma}} = - \lambda \mathcal{P}_{\hspace{-0.08cm}\dot{\gamma}}$ along $\gamma$, where $\mathcal{P}_{\hspace{-0.08cm}\dot{\gamma}}$ is the orthogonal projection onto the orthogonal complement of the velocity field $\dot{\gamma}$. 
\end{lemma}
\textit{Proof.}
Write $\jota(t) =e^{\frac{1}{n-1} \lna \dete \Jxi (t) }$ and differentiate $\mathscr{J}$ to get 

\begin{flalign}
\hspace{0.8cm} \jota^{\prime}
	& = \frac{1}{n-1} \bigg(\frac{d}{dt} \lna \dete \Jxi \bigg) \jota  \notag \\
	& = \frac{1}{n-1} \big( \tra U \big) \jota \notag. &	
\end{flalign}
Hence,

\begin{flalign}
\hspace{0.8cm} \jota^{\prime \prime}
	& = \bigg(-\frac{1}{n-1} \tra U^{2} - \frac{1}{n-1} \tra \textsl{R}_{\dot{\gamma}} \bigg) \jota + \frac{1}{n-1} \big(\tra U\big)  \jota^{\prime} \notag \\
	& \leq - \frac{1}{(n-1)^2} (\tra U)^2 \jota - \frac{1}{n-1} \text{Ric}(\dot{\gamma}, \dot{\gamma})\jota + \frac{1}{(n-1)^2} (\tra U)^2  \jota  \notag \\
	& \leq  (\lambda \circ \distP)(\gamma)\jota. \notag &
\end{flalign}
The last statement follows from the Example \ref{example1}. The initial conditions are clear.$\hspace{0.3cm} \square$ \vspace{0.2cm}

In order to avoid a notational inconvenience, we define the letter $\textit{h}$ to be
\begin{flalign}
\label{def_hit}
\hspace{0.8cm} \hit & \doteq \frac{\has}{n-1}. \tag{$\star$}  &	
\end{flalign}

\begin{lemma}
\label{lemma_comparing_ywithJandj} Let $\textsl{y}$ and $j$ be, respectively, the unique solutions to the following initial value problems (IVP)

\begin{alignat*}{2}
\hspace{-2.5cm}
   & \begin{aligned} & \begin{cases}
  \textsl{y}^{\prime \prime}(t) -\lambda(t) \textsl{y}(t) =0, & \\
 	\textsl{y}(0)=1, \hspace{0.2cm} \textsl{y}^{\prime}(0)=\vert \hit \vert; \hspace{1.5cm} & \\
  \end{cases}\\
  \MoveEqLeft[-1]\text{}
  \end{aligned}
    & \hskip 6em &
  \begin{aligned}
  & \begin{cases}
  j^{\prime \prime}(t) =0, & \\
 	j(0)=1, \hspace{0.3cm} j^{\prime}(0)= \vert \hit \vert. & \\
  \end{cases} \\
  \MoveEqLeft[-1]\text{}
  \end{aligned}
\end{alignat*}
Then, $j \leq \textsl{y}$ on $[0, \infty)$. Moreover, if $\jota$ satisfies the conditions in (\ref{inequalities_satisfies_J}) and $\has = \tra U (0)$, we have \vspace{0.2cm}

\begin{flalign*}
\hspace{0.8cm} \jota &\leq \textsl{y}, \hspace{0.2cm} \text{on } \hspace{0,2cm} [0,t_{\textsl{f}}(\xi)).&
\end{flalign*}

\end{lemma}
\textit{Proof.} First, observe that $\textsl{y}$ is a positive function. In fact, suppose there is a first time $T>0$ for which $\textsl{y}(T)=0$. Then, $\textsl{y}^{\prime \prime}(t)=\lambda(t)\textsl{y}(t)\geq 0$ on $[0,T]$ so that $\textsl{y}^{\prime}(t) = \vert\hit\vert +\int_{0}^{t}\lambda(s)\textsl{y}(s)ds \geq 0$ and $\textsl{y}$ is nondecreasing on $[0, T)$ which leads to a contradiction. This line of reasoning shows that $\textsl{y}^{\prime}(t) = \vert \hit \vert + \int_{0}^{t}\lambda \textsl{y} \geq  \vert \hit \vert$ whence $\textsl{y}(t)-1 \geq t \vert \hit \vert$. This proves that $j \leq \textsl{y}$. To see the second statement, differentiate the quotient $\jota \big/ \textsl{y}$ to get 

\begin{flalign*}
\hspace{0.8cm} \bigg( \frac{\jota(t)}{\textsl{y}(t)} \bigg)^{\prime} &= \frac{\jota^{\prime}(t)\textsl{y}(t)-\jota(t) \textsl{y}^{\prime}(t)}{\textsl{y}(t)^2}\leq 0,&
\end{flalign*}
since $\big( \jota^{\prime}\textsl{y}-\jota \textsl{y}^{\prime} \big)^{\prime}(t) = \jota^{\prime \prime}(t) \textsl{y}(t)  -\jota(t) \textsl{y}^{\prime \prime}(t) \leq 0$ and $(\jota^{\prime}\textsl{y}-\jota \textsl{y}^{\prime})(0)\leq 0$. Thus, $t \mapsto \jota(t)/\textsl{y}(t)$ is nonincreasing and the fact that $(\jota/\textsl{y})(0)=1$ yields $\jota \leq \textsl{y}$.$\hspace{0.3cm} \square$ \vspace{0.2cm}

\begin{com} Here we bring the attention to the function $(\lambda \circ \distP \circ \gamma)$. If the fixed geodesic $\gamma$ is getting far away from the base point $\pzero$ then it may be the case that $(\distP \circ \gamma )$ is contained in a region where the associated function $\lambda$ vanishes. For example, if the function $\lambda$ has compact support and $\Omega$ is a geodesic ball centered at $\pzero$ with radius sufficiently large so that $\distP \circ \gamma (t) \subset \mathbb{R}\setminus \text{supp}(\lambda)$, for all $t$ in the domain of $\gamma$, then all the comparisons in this sections are meaningless. There are other, more general, situations in which this problem arises. We overcome this problem by assuming either that $\lambda(t)\neq 0$, for all $t \in \mathbb{R}^{+}$, or that $\Sigma$ is contained in a region where $(\lambda \circ \distP)(x)$ never vanishes, for all $x \in \Sigma$. The more general situations can be handled as in our proofs with mild modifications. \vspace{0.2cm}
\end{com}

\begin{prop}
\label{proposition_shun_hui_zhu} Let $\textsl{y}$ be the unique solution to the IVP in Lemma \ref{lemma_comparing_ywithJandj}. Then 

\begin{flalign}
\label{inequalities_first_simple:y}
\hspace{0.8cm} \frac{\textsl{y}}{\vert \hit \vert + \int_{0}^{t}\textsl{y} \lambda}& \leq t + \frac{1}{\vert \hit \vert + \int_{0}^{t}\textsl{y} \lambda},&
\end{flalign}
holds only for $t >0$. Multiplying inequality (\ref{inequalities_first_simple:y}) by $\lambda$, integrating on the interval $[\epsilon, t]$, for $0< \epsilon < t$, and using the inequality $\textsl{y} \geq \vert \hit \vert t +1$ from the previous lemma, we get

\begin{flalign}
\label{second_inequality_simple}
\hspace{0.8cm}\lna \bigg( \vert \hit \vert + \int_{0}^{t} \lambda(s) \textsl{y}(s)ds  \bigg) & \leq \int_{\epsilon}^{t}s\lambda (s)ds + \int_{\epsilon}^{t} \frac{\lambda (s)ds}{ \vert \hit \vert + \int_{0}^{s} \lambda(u) (\vert \hit \vert u+1)du} + \lna \bigg( \vert \hit \vert + \int_{0}^{\epsilon}\lambda(s) \textsl{y}(s)ds \bigg) &
\end{flalign}
\end{prop}
\textit{Proof.}
To prove the statements, we proceed as in \cite{10.2307/2374996}. Integrating the equation $\textsl{y}^{\prime \prime} - \lambda \textsl{y} =0$ we get

\begin{flalign}
\label{y_derivative_equation}
\hspace{0.8cm} \textsl{y}^{\prime}(t)& =  \vert \hit \vert + \int_{0}^{t}\lambda(s) \textsl{y}(s) ds, &
\end{flalign}
which holds for all $t \geq 0$ so, integrating once again in $[0,t]$ with $t>0$ and using that $\lambda \neq 0$, we get

\begin{flalign}
\label{value_yt}
\hspace{0.8cm} \textsl{y}(t)-\textsl{y}(0) 
	& = t \vert \hit \vert + \int_{0}^{t}\int_{0}^{u} \lambda(u)\textsl{y}(u) du ds \notag \\
	& = t \vert \hit \vert + \int_{0}^{t}\bigg( \int_{u}^{t} \lambda(u)\textsl{y}(u) ds \bigg) du \notag \\
	& = t \vert \hit \vert  + \int_{0}^{t} (t-u) \lambda(u)\textsl{y}(u)du \\
	& \leq t \vert \hit \vert + t\int_{0}^{t} \lambda(u)\textsl{y}(u)du \notag & 
\end{flalign}
which gives

\begin{flalign*}
\hspace{0.8cm} \frac{\textsl{y}(t)-1}{\vert \hit \vert + \int_{0}^{t} \lambda(s)\textsl{y}(s)ds} & \leq t& 
\end{flalign*}
which is equivalent to (\ref{inequalities_first_simple:y}). Next, multiplying it by $\lambda$ we get

\begin{flalign}
\label{inequality_simple_multiplied_lambda}
\hspace{0.8cm} \frac{\lambda \textsl{y}}{\vert \hit \vert + \int_{0}^{t}\textsl{y} \lambda}& \leq t\lambda  + \frac{\lambda }{\vert \hit \vert + \int_{0}^{t}\textsl{y} \lambda}.&
\end{flalign}
Now, from Lemma \ref{lemma_comparing_ywithJandj}, we have $t\vert \hit \vert +1 \leq \textsl{y}$ so that

\begin{flalign*}
\hspace{0.8cm}
\frac{\lambda (t)}{ \vert \hit \vert + \int_{0}^{t} \lambda(s) \textsl{y}(s)ds}&\leq \frac{\lambda (t)}{ \vert \hit \vert+  \int_{0}^{t} \lambda(s) ( \vert \hit \vert s+1)ds}. &
\end{flalign*}
Thus, integrating (\ref{inequality_simple_multiplied_lambda}) on $[\epsilon ,t]$, for $0 < \epsilon <t$, and using the above inequality gives (\ref{second_inequality_simple}).$\hspace{0.3cm} \square$ \vspace{0.2cm}

\begin{prop}
\label{prop_inequalities_involving_JXy}
Let $\textsl{y}$ be the unique solution to the IVP in Lemma \ref{lemma_comparing_ywithJandj}. Then

\begin{flalign}
\label{inequality_for_y}
\hspace{0.8cm}
\textsl{y}(t) & \leq \big( e^{b_0}ft +1 \big), \hspace{0.2cm} \text{for}\hspace{0.2cm} 0 \leq  t < \infty &
\end{flalign}
where $f = \vert \hit \vert (1+b_0) + b_1$, and $b_0$ and $b_1$ are defined by (\ref{definition_b_0}) and (\ref{definition_b_1}) respectively.
\end{prop}
\textit{Proof.} Let $\epsilon >0$ be a real number (a fixed parameter), and define an auxiliary function $\rho_{\epsilon}:[\epsilon, \infty) \to \mathbb{R}$ by

\begin{flalign*}
\hspace{0.8cm}
\rho_{\epsilon}(t) & =  \int_{\epsilon}^{t} s \lambda(s)ds + 
\int_{\epsilon}^{t} \frac{\lambda(s) ds }{\vert \hit \vert + \int_{0}^{s}\lambda(u) (\vert \hit \vert u+1)du} + \lna \Big( \vert \hit \vert  + \int_{0}^{\epsilon} \lambda(s)\textsl{y}(s)ds \Big). & 
\end{flalign*}

We now distinguish two cases. 

\hspace{0.5cm}\textit{Case 1. ---} Assume first that $ \hit \neq 0$. In this case we may choose the parameter to be $\epsilon =0$. The heart of the matter is that we need to better control the second term on the right hand side of (\ref{second_inequality_simple}). Considering that term, integration by parts show us that  

\begin{flalign*}
\hspace{0.8cm} \int_{0}^{t} \bigg(\frac{1}{1+s \vert \hit \vert} \bigg) \frac{(1+s \vert \hit \vert)\lambda (s)}{ \vert \hit \vert+  \int_{0}^{s} \lambda(u) ( \vert \hit \vert u+1)du} ds \notag \vspace{0.8cm} 
	&  = \frac{1}{1+t \vert \hit \vert} \lna\bigg( \vert \hit \vert +\int_{0}^{t} \lambda(s) ( \vert \hit \vert s+1)ds \bigg) \notag \\
	&  -\lna \vert \hit \vert - \int_{0}^{t} \frac{- \vert \hit \vert }{(1+s \vert \hit \vert)^{2}} \lna\bigg(  \vert \hit \vert +\int_{0}^{s} \lambda(u) ( \vert \hit \vert u+1)du \bigg) ds \notag \vspace{0.8cm} \\
	& \leq  \frac{1}{1+t \vert \hit \vert} \lna\bigg(  \vert \hit \vert + \int_{0}^{t} \lambda(s) ( \vert \hit \vert s+1)ds \bigg) \notag \\
	& -\lna  \vert \hit \vert + \vert \hit \vert \lna f\int_{0}^{t} \frac{1}{(1+s \vert \hit \vert)^{2}} ds. \notag & \\
\end{flalign*}
Observe that $\rho_{0}^{\prime}(t) \geq 0$ so that $\rho_{0}$ is nondecreasing and note that $\rho_{0}$ is bounded above by $b_0 + \lna f$. Hence, by (\ref{y_derivative_equation}) and analysing inequality (\ref{second_inequality_simple}) with $\epsilon = 0$ and $\hit \neq 0$ we conclude that

\begin{flalign*}
\hspace{0.8cm}\textsl{y}^{\prime}(t)  & \leq e^{ \rho_{0}(t)} \leq e^{b_0} f &
\end{flalign*}
holds for all $t \geq 0.$

\hspace{0.5cm}
\textit{Case 2. ---} We now suppose that $ \hit =0$. In this case, the function $\rho_{\epsilon}$ becomes

{\small
\begin{flalign}
\label{rho_eps_2_cas}
\hspace{0.8cm}\rho_{\epsilon}(t) & = \int_{0}^{t}s\lambda (s)ds + \lna \bigg( \int_{0}^{t}\lambda(s)ds \bigg) - \lna \bigg( \int_{0}^{\epsilon} \lambda (s)ds \bigg)  +  \lna  \bigg(\int_{0}^{\epsilon} \lambda(s) \textsl{y}(s)ds \bigg) -  \int_{0}^{\epsilon} s\lambda(s)ds. & 
\end{flalign}}
Again, $\rho_{\epsilon}^{\prime} \geq 0$ and $\rho_{\epsilon}$ is nondecreasing. In addition, letting $t \rightarrow \infty$ we see from (\ref{rho_eps_2_cas}) that $\rho_{\epsilon}$ satisfies 

\begin{flalign*}
\hspace{0.8cm}
\rho_{\epsilon}(t) & \leq b_0 + \lna b_1 + \lna\bigg( \frac{\int_{0}^{\epsilon} \lambda (s) \textsl{y}(s)ds}{\int_{0}^{\epsilon}\lambda (s)ds }\bigg)- \int_{0}^{\epsilon}s\lambda(s)ds. & 
\end{flalign*}
Therefore, by (\ref{y_derivative_equation}), after applying the exponential function, we may rewrite inequality (\ref{second_inequality_simple}) as

\begin{flalign}
\label{inequality8_one_of_main}
\hspace{0.8cm}
\textsl{y}^{\prime}(t) & \leq  e^{\rho_{\epsilon}(t)} \leq \bigg( e^{b_0}b_1  \frac{1}{e^{ \int_{0}^{\epsilon}s\lambda(s)ds}}  \bigg)
\frac{\displaystyle{\int_{0}^{\epsilon} \lambda (s) \textsl{y}(s)ds}}{\displaystyle{ \int_{0}^{\epsilon}\lambda(s)ds }} , &
\end{flalign}
which holds only for $t>0$. Since the LHS of (\ref{inequality8_one_of_main}) does not depend on $\epsilon$ and it holds for all $\epsilon >0$ small, we may take the limit on the RHS, if it exists, when $\epsilon \rightarrow 0$. But, the quotient of integrals above clearly satisfies 

\begin{flalign*}
\hspace{0.8cm}
\lim_{\epsilon \rightarrow 0^{+}} \frac{1}{\displaystyle{ \int_{0}^{\epsilon}\lambda }} \int_{0}^{\epsilon} \lambda \textsl{y}& = \lim_{\epsilon \rightarrow 0^{+}} \textsl{y}(\epsilon)= 1. &
\end{flalign*}
Thus, letting $\epsilon \rightarrow 0$ in (\ref{inequality8_one_of_main}) gives $\textsl{y}^{\prime} \leq e^{b_0}b_1$. As a result, regardless of the signal of $\hit$ we have $\textsl{y}^{\prime} \leq e^{b_0}f$. Integrating this differential inequality gives (\ref{inequality_for_y}).$\hspace{0.3cm} \square$ \vspace{0.2cm}

\begin{com}
\label{dec_asy_beh}
In view of Proposition \ref{prop_inequalities_involving_JXy}, it is natural to wonder about the asymptotic behavior of $t \mapsto \textsl{y}(t)/(e^{b_0}ft+1)$. As we will see later, it will play a crucial role in understanding the connection between the decay of $\lambda$ and the rigidity statement. To understand this issue, we introduce the notation $\mathbb{X}(t) =  (e^{b_0}ft+1)$. By equations (\ref{value_yt}) and (\ref{y_derivative_equation}) we have

\begin{flalign*}
\hspace{0.8cm}
	\bigg(\frac{\textsl{y}(t)}{\mathbb{X}(t)} \bigg)^{\prime}
	& = \frac{e^{b_0}f}{\mathbb{X}(t)^2}\bigg[ t\bigg(\vert \hit \vert + \int_{0}^{t}\lambda(s) \textsl{y}(s) ds\bigg) - \textsl{y}(t) \bigg] + \frac{1}{\mathbb{X}(t)^2} \bigg[\vert \hit \vert + \int_{0}^{t}\lambda(s) \textsl{y}(s) ds \bigg] \\
		& = \frac{e^{b_0}f}{\mathbb{X}(t)^2}\bigg[ -1 + \int_{0}^{t}u\lambda(u) \textsl{y}(u) du \bigg] + \frac{1}{\mathbb{X}(t)^2} \bigg[\vert \hit \vert + \int_{0}^{t}\lambda(s) \textsl{y}(s) ds \bigg].&
\end{flalign*}
So, on the one hand, the signal of the above derivative depends on the factor $-1 + \int_{0}^{t}u\lambda(u) \textsl{y}(u) du$. On the other hand, it depends on the numerator $\textsl{y}^{\prime}\mathbb{X}-\mathbb{X}^{\prime}\textsl{y}$ which is increasing as long as $\lambda$ does not vanish since

\begin{flalign*}
\hspace{0.8cm}
\big(\textsl{y}^{\prime}(t)\mathbb{X}(t)-\mathbb{X}^{\prime}(t)\textsl{y}(t) \big)^{\prime} & = \lambda(t) \textsl{y}(t) \mathbb{X}(t) \geq 0. &
\end{flalign*}
In particular,

\begin{flalign}
\label{thet_beg_dec}
\hspace{0.8cm}
\frac{d}{dt}\bigg\vert_{t=0} \bigg(\frac{\textsl{y}(t)}{\mathbb{X}(t)} \bigg)& = \vert h \vert -e^{b_0}f <0.&
\end{flalign}
\end{com}

\section{A sharp geometric inequality, \textit{\nth{1}} case}
\label{sharp_dom_asym_non}

In this section we give the proof of the Theorem \ref{main_theo}. First we slightly deviate from the course of the proof to establish some notations and terminology. \vspace{0.2cm}

\begin{dig}
\label{dig_on_collective_behavior} We now look at the collective behavior of geodesics that start at the boundary $\Sigma =\partial \Omega$ and realize the distance $r(x)= \text{dist}_{\textsl{g}}(x, \Sigma)$, for $x \in \textsl{M} \setminus \Omega$ to gather all the information developed so far into single maps. For each $x \in \Sigma$, define

\begin{flalign}
\label{definition_of_X_f}
\hspace{0.8cm} \mathbb{X}_x(t) & = e^{b_0} f(x)t +1, \hspace{0.2cm} t \in \mathbb{R},&
\end{flalign}
where $f: \Sigma \to \mathbb{R}$ is the function $f(x) = \vert \hit \vert(x)(1+b_0) + b_1$ and $b_0$ and $b_1$ are defined by (\ref{definition_b_0}) and (\ref{definition_b_1}) respectively. Remind that $\hit(x) = \has(x)/(n-1)$ as defined in (\ref{def_hit}), where $\has(x)$ is the mean curvature of $\Sigma$ at $x$. Recall that 

\begin{flalign*}
\hspace{0.8cm}\jota_x(t)& = \dete J_{\xi(x)}(t)^{\frac{1}{n-1}}, \hspace{0.2cm} t \in [0, \txibx),&
\end{flalign*}
as defined in (\ref{definition_of_J}), where $\xi(x) = \xi_x$ is the outward, unit and normal vector at $x \in \Sigma$ and $\txibx$ is the first time the geodesic $\gamma_{\xi_x}$ hits a focal point. We sometimes consider parameters as variables and write, e.g., $\mathbb{X}(x,t)$ instead of $\mathbb{X}_x(t)$ or $\jota(x,t)$ in place of $\jota_{x}(t)$. By Proposition \ref{prop_inequalities_involving_JXy} and Lemma \ref{lemma_comparing_ywithJandj} we have $\jota_{x}^{n-1}(t) \leq \mathbb{X}_{x}^{n-1}(t)$, $t \in [0, \txibx)$. Next, define the function $\hat{\theta}  :\Sigma \times [0, \infty) \to \mathbb{R}$ by 

\begin{flalign}
\label{definition_of_theta_x}
\hspace{0.8cm} \hat{\theta} (x,t) & = \bigg(\frac{\jota_{x}(t)}{\mathbb{X}_{x}(t)}\bigg)^{n-1}. &
\end{flalign}
For each $x \in \Sigma$ the quantity $\hat{\theta} (x,\cdot )$ is clearly positive on $[0, \txibx)$. Also, if $\textsl{y}_{x}$ denotes the unique solution to $\textsl{y}_{x}^{\prime \prime}(t)-\lambda(t)\textsl{y}_{x}(t)=0$, with suitable initial conditions, where $\lambda(t)$ is implicitly understood to be the function $\lambda(\distP(\gamma_{\xi_x}(t)))$, then by Proposition \ref{prop_inequalities_involving_JXy} we have

\begin{flalign}
\label{fundamental}
\hspace{0.8cm} \hat{\theta}(x, t)^{\frac{1}{n-1}}& \leq \frac{\jota_{x}(t)}{\textsl{y}_x(t)}, \hspace{0.2cm} x \in \Sigma,\hspace{0.2cm}  t \in  [0, \txibx),&
\end{flalign}
where $t \mapsto \jota_x(t)/\textsl{y}_x(t)$ is monotone decreasing by Lemma \ref{lemma_comparing_ywithJandj}. Moreover, for any $0 \leq a < b \leq \infty$ such that $b \leq \tau_c(x)$, where $\tau_c(x)$ is the cut time defined in (\ref{cut_time}), we shall consider the supremum of $\hat{\theta} (x,\cdot)$ for $t \in [a,b)$. That supremum is denoted by

\begin{flalign*}
\hspace{0.8cm} \sup_{[a, b)} \hat{\theta} (x) \doteq & \sup_{t \in [a, b)} \hat{\theta} (x, t) \leq 1, \hspace{0.2cm} x \in \Sigma . &
\end{flalign*}
Finally we define the function whose asymptotic behavior is key for this paper. Let $\theta :\Sigma \times [0, \infty) \to \mathbb{R}$ be defined by

\begin{flalign*}
%\label{}
\hspace{0.8cm} \theta(x,t)^{\frac{1}{n-1}} & =\frac{\textsl{y}_{x}(t)}{\mathbb{X}_{x}(t)}. &
\end{flalign*}
In Remark \ref{dec_asy_beh} we analysed the first derivative of $\theta_x(t)=\theta(x,t)$. In particular, if $\lambda \neq 0$ then by (\ref{thet_beg_dec}) we have

\begin{flalign}
\label{fir_der_zero_tht}
\hspace{0.8cm}
\frac{\1}{\1 t }\bigg\vert_{t=0} & \theta(x,t)^{\frac{1}{n-1}} < 0,&
\end{flalign}
so that for all $x \in \Sigma$, the function $t \mapsto \theta_x(t)^{\frac{1}{n-1}}$ is decreasing on a neighborhood of zero. $\hspace{0.3cm} \square$
\end{dig}

Next, we proceed with the proof of Theorem \ref{main_theo}, which is divided into two steps. To prove inequality (\ref{willmorex_type_inequality_main}) we may assume, without loss of generality, that $\Omega$ has no hole, i.e., $\textsl{M}\setminus \Omega$ has no bounded component. \vspace{0.2cm}

\textbf{Proof of inequality (\ref{willmorex_type_inequality_main}).}
\label{prf_ineq_four} 
Suppose first that $\lambda = 0$. Then, by hypothesis, $(\textsl{M},\textsl{g})$ has nonnegative Ricci curvature and, therefore, the argument in \cite{AFST_2023_6_32_1_173_0} goes through without modification to conclude inequality (\ref{willmorex_type_inequality_main}). We now assume that $\lambda \neq 0$. Thus, all the comparison arguments in Section \ref{section_elel_ineq} hold. \vspace{0.1cm}

Denote by $\mathcal{T}_{\Omega}(R)$ the geodesic tube of radius $R$ about $\Omega$. Then, for any $R>0$, we have

\begin{flalign*}
\hspace{0.8cm} \textit{vol}(\mathcal{T}_{\Omega}(R))
	& = \vert \Omega \vert + \int_{\Sigma} \int_{0}^{R \wedge \tau_c(x)} \dete J_{\xi(x)} (t) dt d\sigma (x) \notag \\ 
	& \leq \vert \Omega \vert + \int_{\Sigma} \int_{0}^{R \wedge \tau_c(x)} \mathbb{X}(x,t)^{n-1} dt d\sigma (x)  \notag \\
	& \leq \vert \Omega \vert + \int_{\Sigma} \int_{0}^{R} \big( e^{b_0} f(x)t +1 \big)^{n-1} dt d\sigma(x) \notag \\
	& = \vert \Omega \vert + e^{(n-1)b_0} \frac{R^{n}}{n} \int_{\Sigma} f(x)^{n-1}   d\sigma(x) +\mathcal{O}(R^{n-1}),&
\end{flalign*}
where $R \wedge \tau_c$ denotes the minimum between $R$ and $\tau_c$. Dividing both sides of the inequality above by $\wn R^{n} = \vert \mathbb{S}^{n-1} \vert R^{n}/n$ and letting $R \rightarrow \infty$ we obtain

\begin{flalign*}
\hspace{0.8cm} AVR(\textsl{g}) & \leq e^{(n-1)b_0} \frac{1}{\vert \mathbb{S}^{n-1} \vert} \int_{\Sigma} f(x)^{n-1}  d\sigma(x)  &
\end{flalign*}
which implies (\ref{willmorex_type_inequality_main}). 
\vspace{0.2cm}

\textbf{Equality discussion.} We first prove the case \textit{($\mathbf{\textsl{W}}$\hspace{-0.02cm}\textbf{1})}. If $\lambda = 0$ then $(\textsl{M}, \textsl{g})$ has nonnegative Ricci curvature. If equality in (\ref{willmorex_type_inequality_main}) holds, we have

\begin{flalign*}
\hspace{0.8cm}
\int_{\Sigma}\bigg \vert \frac{\has}{n-1} \bigg \vert^{n-1} d\sigma &= AVR(\textsl{g}) \vert\mathbb{S}^{n-1} \vert.&
\end{flalign*}
Then the argument in \cite{AFST_2023_6_32_1_173_0} goes through \textit{ipsis litteris} to conclude this first part. Next, we move to the case \textit{($\mathbf{\textsl{W}}$\hspace{-0.02cm}\textbf{2})}. Suppose we have

\begin{flalign}
\label{equality_discussion_asymptotic_curvature}
\hspace{0.8cm}
e^{(n-1)b_0}\int_{\Sigma}\bigg( \bigg \vert \frac{ \has}{n-1} \bigg \vert (1+b_0) + b_1  \bigg)^{n-1} d\sigma &= AVR(\textsl{g}) \vert\mathbb{S}^{n-1} \vert.&
\end{flalign}
We must show that $\Sigma$ is a totally umbilical hypersurface with constant mean curvature on components. Observe first that from the continuity of $\tau_c$ and from the inequalities in \hyperref[prf_ineq_four]{Proof of inequality (4)} above we deduce $\tau_c \equiv \infty $ on $\Sigma$. Now, on the one hand, for all $x \in \Sigma$ and $R^{\prime}>0$, the hypothesis $\lambda \neq 0$ together with (\ref{fundamental}) ensures that

\begin{flalign}
\label{chn_hat_tht_vs}
\hspace{0.8cm}
\sup_{[R^{\prime}, \infty)} \hat{\theta} (x)=
\sup_{[R^{\prime}, \infty)} \bigg(\frac{\jota_{x}(t)}{\mathbb{X}_{x}(t)}\bigg)^{n-1} & \leq \sup_{ [R^{\prime}, \infty)} \bigg(\frac{\jota_{x}(t)}{\textsl{y}_{x}(t)}\bigg)^{n-1} = \bigg(\frac{\jota_{x}(R^{\prime})}{\textsl{y}_{x}(R^{\prime})}\bigg)^{n-1} \leq 1. &
\end{flalign}
On the other hand, for any $0 < R^{\prime} < R$, we have

\begin{flalign*}
%\label{string1}
\hspace{0.8cm} \textit{vol}(\mathcal{T}_{\Omega}(R))
	& = \vert \Omega \vert + \int_{\Sigma} \int_{0}^{R} \dete J_{\xi(x)}(t) dt d\sigma (x)  \\
	& = \vert \Omega \vert + \int_{\Sigma} \int_{R^{\prime}}^{R}  \hat{\theta} (x,t)\mathbb{X}(x,t)^{n-1}dt d\sigma (x) + 
\int_{\Sigma} \int_{0}^{R^{\prime}}\dete J_{\xi(x)}(t) dt d\sigma (x) \notag \\
	& \leq  
\vert \Omega \vert + \int_{\Sigma} \sup_{[R^{\prime}, R]} \hat{\theta} (x)  \int_{R^{\prime}}^{R}\big(e^{b_0} f(x)t +1\big)^{n-1}dt d\sigma (x) + \int_{\Sigma} \int_{0}^{R^{\prime}}\dete J_{\xi(x)}(t) dt d\sigma (x) \notag \\
	& = \vert \Omega \vert + e^{(n-1)b_0} \frac{R^{n}}{n} \int_{\Sigma} \sup_{[R^{\prime}, R]} \hat{\theta} (x) f(x)^{n-1}   d\sigma(x) + \mathcal{O}(R^{n-1}). \notag &
\end{flalign*}
Dividing both sides by $\vert \mathbb{B}^{n}\vert  R^{n} \big/ \vert \mathbb{S}^{n-1} \vert = R^{n}/n$ and letting $R\rightarrow \infty$, we get

\begin{flalign}
\label{avr_againsttheta}
\hspace{0.8cm}
 \vert\mathbb{S}^{n-1} \vert AVR(\textsl{g}) \leq  & e^{(n-1)b_0}\int_{\Sigma}\sup_{[R^{\prime}, \infty)} \hat{\theta}(x)\bigg( \bigg \vert \frac{ \has(x)}{n-1} \bigg \vert (1+b_0) + b_1  \bigg)^{n-1}d\sigma(x), &
\end{flalign}
for all $R^{\prime}>0$. Subtracting (\ref{equality_discussion_asymptotic_curvature}) from (\ref{avr_againsttheta}), one easily deduces that for almost every $x \in \Sigma$ and $R^{\prime}>0$

\begin{flalign*}
\hspace{0.8cm}
\sup_{[R^{\prime}, \infty)} \hat{\theta}(x) - 1 & = 0. &
\end{flalign*}
Hence, according to (\ref{chn_hat_tht_vs}), we have $\jota_{x}(R^{\prime})\big/ \textsl{y}_{x}(R^{\prime})=1$ for almost every $x \in \Sigma$ and for all $R^{\prime}>0$. By continuity, this holds everywhere. In particular, $\hat{\theta} \equiv \theta$. The initial conditions on $\jota_{x}(t)$ and $\textsl{y}_{x}(t)$ return $\vert \has \vert = \has \geq 0$. Inspecting the comparison argument above and reasoning as in Lemma (\ref{lemma1_elel}), we get that

\begin{enumerate}[label=\textit{(E\hspace{-0.035cm}\arabic*)},leftmargin=1.7cm]
\item $\mathcal{H}_{r} = \frac{\Delta r}{n-1} \mathcal{P}_{\hspace{-0.08cm}\1_r}$;

\item $\Ricc(\1_r,\1_r)=-(n-1)(\lambda \circ \distP)$,
\end{enumerate}
holds along $\Phi([0, \infty) \times \Sigma)$, where $\distP$ in $\textit{(E\hspace{-0.035cm}2)}$ above is evaluated at the integral curves of $\grad r$. It follows from the first equation $\textit{(E\hspace{-0.035cm}1)}$ that $\Sigma$ is umbilic. It follows from the second equation that $\xi_x = \1_r\big\vert_x$ must be an eigenvector of the Ricci tensor of $(\textsl{M}, \textsl{g})$ since $\Ricc\geq -(n-1)(\lambda \circ \distP)\textsl{g}$, for $x \in \Sigma$. Therefore $\Ricc(\xi_x,e_k)=0$, where $\{ e_{1},..., e_{n-1}\}$ is an orthonormal basis for the tangent space $T_x\Sigma$. Using the Codazzi equations, we deduce that the curvature tensor of $(\textsl{M}, \textsl{g})$ satisfies $\textsl{R}(e_i, e_j, e_k, \xi) = \frac{1}{(n-1)}\big(e_i \has \delta_{jk}  - e_j \has \delta_{ik}\big)$ so that 

\begin{flalign*}
\hspace{0.8cm}
0 =&  \Ricc(e_j, \xi) = \frac{(n-2)}{n-1}e_j\has. &
\end{flalign*}
Consequently, $\has$ is locally constant. This establishes the first claim. Next, Suppose $\Sigma$ is connected. Write $\mathcal{V} = \{ t\xi_x: x\in \Sigma, \hspace{0.1cm} t \in [0, \infty) \}$ and for each $t \in [0, \infty)$, define $\Sigma_{t}^{*} = r^{-1}(t) \cap \expo^{\perp}(\mathcal{V})$, which is the part of the level sets consisting of regular points, and because  $\tau_c \equiv \infty$, it is equal to $r^{-1}(t)$ itself. Moreover, since $r^{-1}[a,b]$ is compact and contains no critical points of $r$ we conclude, by Morse theory, that the ${\Sigma_{t}^{*}}'\textsl{s}$ are connected, whence $\textsl{M}$ has only one end. \vspace{0.01cm}

For each $t \in [0, \infty)$, let $\has[t]$ be the mean curvature of $\Sigma_{t}^{*}$ and let $\varphi$ be the unique solution to $\varphi^{\prime}(t)- \frac{\has[t]}{n-1}\varphi(t)=0$, with initial condition $\varphi(0)=1$, which is well posed as long as $r^{-1}(t)$ is connected. It is easy to see that $\varphi$ is increasing and $\therefore$ positive (cf. Lemma \ref{lemma_comparing_ywithJandj}). Now, let us show that $\varphi$ satisfies a Jacobi type equation and that $\Sigma$ is a geodesic sphere centered at $\pzero$. To obtain these claims, observe that the function $r \mapsto \frac{\varphi^{\prime}(r)}{\varphi(r)}$ satisfies a Riccati type equation since, by the Bochner formula, we have
\begin{flalign}
\label{uniq_sols}
\hspace{0.8cm}
\frac{1}{2} \Delta \vert \nabla r \vert^{2}
	& = \vert D^{2}r \vert^2 + \langle \grad r, \grad \Delta r \rangle + \Ricc (\grad r, \grad r) \notag \\
	& = \frac{(\Delta r)^{2}}{n-1} + \frac{\1}{\1 r} \Delta r - (n-1)(\lambda \circ \distP)(\gamma_{\xi_x}(r)) \notag \\
	& = (n-1)\bigg(  \Big(\frac{\varphi^{\prime}}{\varphi}\Big)^{\prime} + \Big(\frac{\varphi^{\prime}}{\varphi}\Big)^{2} - (\lambda \circ \distP)\big(\gamma_{\xi_x}(r)\big)\bigg)=0,&
\end{flalign}
for all $x\in \Sigma$ and each integral curve $\gamma_{\xi_x}$ of $\grad r$. By uniqueness of solutions to (\ref{uniq_sols}) with given initial data, the function $\Sigma \ni x \mapsto \distP(x)$ can only be constant. Therefore, $\Omega$ is a geodesic ball centered at $\pzero$ with radius $r_0 = \text{dist}_{\textsl{g}}(\pzero, \Sigma)$. Moreover, $\textsl{M} \setminus \Omega$ is easily seen to be locally a warped product, since defining $u: \textsl{M} \setminus \Omega \to \mathbb{R}$ by

\begin{flalign*}
\hspace{0.8cm}
u(w)& \doteq \int_{0}^{r(w)}\varphi(t)dt,  &
\end{flalign*}
gives $\grad u = \varphi(r) \grad r$ so that $D^{2}u = \varphi^{\prime}\textsl{g}$ which fulfils a characterization of warped product with 1-dimensional factor (see e.g.  \cite{cheeger1996lower}, pp. 192 \textit{--} 194). Next, we globalize this result.

\hspace{0.5cm}
\textit{\textbf{Claim.}---} \textit{The manifold $\textsl{M} \setminus \Omega$ is isometric to $\big([r_0, \infty) \times \Sigma , dr\otimes dr + \varrho(r)^2 \textsl{g}_{\Sigma}\big)$}. \vspace{0.05cm}

In fact, since $\Phi$ is a diffeomorphism from $[0, \infty) \times \Sigma$ onto $\textsl{M} \setminus \Omega$, the pulled back metric takes the form

\begin{flalign*}
\hspace{0.8cm}
\Phi^{*}\textsl{g}& = dr \otimes dr + \textit{g}_{r}, &
\end{flalign*}
where $r \mapsto \textit{g}_{r}$ is a $r$-dependent family of metrics over $\Sigma$ with $\textit{g}_{0} = \textsl{g}_{\Sigma}$. By equation $\textit{(E\hspace{-0.035cm}1)}$ again, using appropriate coordinates $\{x^{1},..., x^{n-1}\}$ on $\Sigma$ we may deduce that 

\begin{flalign*}
\hspace{0.8cm}
\frac{1}{2} \frac{\1}{\1 r} (\textit{g}_{r})_{kl} & = \frac{\Delta r}{n-1}(\textit{g}_{r})_{kl} = \frac{\varphi^{\prime}}{\varphi}(\textit{g}_{r})_{kl}  &
\end{flalign*}
Therefore, $\textit{g}_{r} = \varphi(r)^2 \textsl{g}_{\Sigma}$. Applying the translation $\varrho(r) = \varphi (r-r_0)$, $r \geq r_0$ it is straightforward to check that $\varrho$ satisfies the required conditions. In particular, the identity (\ref{iso_asy_to}) is a consequence of Proposition \ref{prp_decay} together with the fact that $\sup\{ \hat{\theta}(x,t): t \in [R^{\prime}, \infty)\}=1$.$\hspace{0.3cm} \square$ \vspace{0.2cm}

\begin{com}
\label{obs_sup_varphi}
Observe that $\Sigma$ being a geodesic sphere has no role in the proof above. Thus, it is reasonable to contest this fact in view of the following reason. A priori, it can be the case that $\lambda$ has compact support and $\{ \distP(x): x \in \Sigma \} \subset \mathbb{R}^{+} \setminus \text{supp}(\lambda)$ so that $\Sigma \ni x \mapsto \lambda(\distP(x)) \equiv 0$. In this case, neither the solution $\varphi$ to the equation $\varphi^{\prime \prime}(t)=0$ nor its translation $\varrho$ satisfy the limit in (\ref{iso_asy_to}). We conclude that the set $\{ \distP(x): x \in \Sigma \}$ must intersect the interior of $\text{supp}(\lambda)$.
\end{com}

\subsection{The decay of $\lambda$}
\label{decay}
Here we address the issue of the decay of the associated function $\lambda$, especially regarding the equality case in the main theorem and in relation with the quantities defined in Digression \ref{dig_on_collective_behavior}.  \vspace{0.2cm}

\begin{prop}
\label{prp_decay}
Suppose we have equality in (\ref{willmorex_type_inequality_main}). If $\lambda \neq 0$ then for every $x_0 \in \Sigma$ it does not exist an open neighborhood $\mathcal{W}$ of $x_0$ in $\Sigma$ such that $t \mapsto \theta(x, t)^{\frac{1}{n-1}}$ does not have critical points for all $x \in \mathcal{W}$. In particular, if for some $x_0 \in \Sigma$, we have

\begin{enumerate}[label=($\mathbf{\textsl{D}}$\hspace{-0.02cm}\textbf{.\arabic*}), leftmargin=1.3cm]
\item $\frac{\1}{\1 t}\theta(x_0, t)^{\frac{1}{n-1}} = 0$ for all $t \geq \tau$ for some $\tau>0$ then $\lambda$ has compact support;
\item $\frac{\1}{\1 t} \theta(x, t)^{\frac{1}{n-1}}>0$, for all $t$ big enough, then $ {\displaystyle \lim_{r \rightarrow \infty}} \frac{\1}{\1 t}\theta(x, r)^{\frac{1}{n-1}}  =0$. 
\end{enumerate}
\end{prop}
\textit{Proof.} Assume we have equality in the inequality (\ref{willmorex_type_inequality_main}). Suppose for the sake of contradiction that there exists $x_0 \in \Sigma$ and some open neighborhood $\mathcal{W}$ of $x_0$ in $\Sigma$ such that $t \mapsto \theta(x, t)^{\frac{1}{n-1}}$ does not have critical points for all $x \in \mathcal{W}$. From (\ref{fir_der_zero_tht}) and the hypothesis, we have $\frac{\1}{\1 t} \theta (x, t)^{\frac{1}{n-1}}<0$ for all $t \geq 0$ so that $t \mapsto \theta(x, t)^{\frac{1}{n-1}}$ is (strictly) decreasing along $[0, \infty)$ and therefore so is $t\mapsto \hat{\theta}(x,t)$, for all $x \in \mathcal{W}$, since we have seen that equality in (\ref{willmorex_type_inequality_main}) implies $\hat{\theta} = \theta$. Hence, for $0 < R^{\prime}< R$, a simple analysis of the asymptotic behavior of $\textit{vol}(\mathcal{T}_{\Omega}(R))/ \vert \mathbb{B}^{n}\vert  R^{n}$, where the numerator satisfies a chain of inequalities of the type

\begin{flalign*}
\hspace{0.8cm}
\textit{vol}(\mathcal{T}_{\Omega}(R))
	& = \vert \Omega \vert + \int_{\mathcal{W}} \int_{0}^{R} \dete J_{\xi(x)}(t)dt d\sigma (x) + \int_{\Sigma \setminus \mathcal{W}} \int_{0}^{R} \dete J_{\xi(x)}(t)dt d\sigma (x)    \\
	& =  \vert \Omega \vert + \int_{\mathcal{W}} \int_{R^{\prime}}^{R} \theta(x,t)X(x,t)^{n-1}dt d\sigma (x) + \int_{\mathcal{W}} \int_{0}^{R^{\prime}} \dete J_{\xi(x)}(t)dt d\sigma (x) + \\ 
	& \hspace{1cm} \int_{\Sigma \setminus \mathcal{W}} \int_{0}^{R} \dete J_{\xi(x)}(t)dt d\sigma (x)    \\
	& < \vert \Omega \vert + \int_{\mathcal{W}}\theta(x,R^{\prime}) \int_{R^{\prime}}^{R} X(x,t)^{n-1}dt d\sigma (x) + \int_{\mathcal{W}} \int_{0}^{R^{\prime}} \dete J_{\xi(x)}(t)dt d\sigma (x) + \\ 
	& \hspace{1cm} \int_{\Sigma \setminus \mathcal{W}} \int_{0}^{R} \dete J_{\xi(x)}(t)dt d\sigma (x)   & \\
\end{flalign*}
yields

\begin{flalign*}
\hspace{0.8cm}
AVR(\textsl{g}) < &
\frac{e^{(n-1)b_0}}{\vert\mathbb{S}^{n-1} \vert} \int_{\Sigma}\bigg( \bigg \vert \frac{ \has}{n-1} \bigg \vert (1+b_0) + b_1  \bigg)^{n-1} d\sigma ,&
\end{flalign*}
contradicting equality in (\ref{willmorex_type_inequality_main}). This proves the first assertion of the proposition. To prove \textit{($\mathbf{\textsl{D}}$\hspace{-0.02cm}\textbf{.1})} suppose $\frac{\1}{\1 t} \theta (x_0, r)^{\frac{1}{n-1}}= 0$ for all $r \geq \tau$ for some positive $\tau$ and some $x_0 \in \Sigma$. Then, $\theta(x_0, t)^{\frac{1}{n-1}} = \textsl{y}_{x_0}(t)/\mathbb{X}_{x_0}(t) = \text{c} \in \mathbb{R}$ for $t \geq \tau$ (it is easy to see that $\text{c}=1$, but this has no role here). Then, we must have $\textsl{y}_{x_0}(t) = \text{c}\mathbb{X}_{x_0}(t)$, for all $t \geq \tau$. Therefore, 

\begin{flalign*}
\hspace{0.8cm}
0 & = \textsl{y}_{x_0}^{\prime \prime}(t) = \lambda(t) \textsl{y}_{x_0}(t) &
\end{flalign*}
for all $t \geq \tau$ and because $\textsl{y}_{x_0} >0$ this proves that $\lambda$ has compact support. To see \textit{($\mathbf{\textsl{D}}$\hspace{-0.02cm}\textbf{.2})} notice that $t \mapsto \theta (x_0, t)^{\frac{1}{n-1}}$ is eventually increasing, by hypothesis. Since $\theta = \hat{\theta}$ and $\sup\{ \hat{\theta}(x,t): t \in [R^{\prime}, \infty)\}=1$, we conclude that the limit of $\theta(x,r)$ for $r \rightarrow \infty$ is equal to $1$ which completes this case. $\hspace{0.3cm} \square$ \vspace{0.2cm}

It is relevant to point out that the first assertion of Proposition \ref{prp_decay} is equivalent to: 
\vspace{0.1cm}

\begin{quote}
\textit{Suppose that $\lambda \neq 0$. If there exists $x_0 \in \Sigma$ and an open neighborhood $\mathcal{W}$ of $x_0$ in $\Sigma$ such that $t \mapsto \theta(x, t)^{\frac{1}{n-1}}$ does not have critical points on $\mathcal{W}$ then the inequality (\ref{willmorex_type_inequality_main}) is strict.}
\end{quote}

\section{A sharp geometric inequality, \textit{\nth{2}} case}
\label{new}

In this section we sketch the proof of Theorem \ref{sec_main_theo}. We first digress about its  preliminary conditions. \vspace{0.2cm}

\begin{dig} $\big[$\textit{The notion of asymptotic$^{\text{\ding{93}}}$ volume ratio} $\big]$ \newline
\label{dig_hom_hyp}
Assume $(\textsl{M},\textsl{g})$ has asymptotically$^{\text{\ding{93}}}$ nonnegative Ricci curvature relatively to $\gamze$, where $\gamze$ is a connected and closed hypersurface. Assume additionally that $\gamze$ separates the ambient $\textsl{M}$, that is, $\textsl{M} \setminus \gamze = \textsl{M}^{+} \cup \textsl{M}^{-}$ where both $\textsl{M}^{+}$ and $\textsl{M}^{-}$ are connected subsets of $\textsl{M}$, at least one of these components is unbounded and the other is possibly empty. Assume $\textsl{M}^{+}$ to be the unbounded component. Take a unit and normal vector field $\nu: \gamze \to \nu \gamze$ that points into $\textsl{M}^{+}$. Our main interest relies on the analysis of the asymptotic behavior of

\begin{flalign*}
\hspace{0.8cm}
\textit{vol}(\mathcal{T}_{\gamze}^{+}(R)) &= \int_{\gamze} \int_{0}^{R \wedge \tau_c(x)} \dete J_{\nu(x)} (t) dt d\sigma_{\hspace{-0.04cm}\textbf{\tiny{0}}} (x),&
\end{flalign*}
where $d\sigma_{\hspace{-0.04cm}\textbf{\tiny{0}}}$ is the induced volume element on $\gamze$. We observe now that all comparisons from Section \ref{section_elel_ineq} go through if we replace the distance $\distP$ by $\text{dist}_{\gamze}$. Indeed, one defines the Jacobian $\mathscr{J}$ just as in formula (\ref{definition_of_J}) and notices that the only difference between the comparisons there and the case at hand is that we are now considering geodesics $\gamma_{\nu_x}$ in $\textsl{M}$ with $\dot{\gamma}_{\nu_x}(0) = \nu_x$ which implies:

\begin{flalign*}
\hspace{0.8cm}
\text{dist}_{\gamze}(\gamma_{\nu_x}(t)) & = t \hspace{0.1cm} \Rightarrow \hspace{0.1cm} \lambda(\text{dist}_{\gamze}(\gamma_{\nu_x}(t))) = \lambda(t), \hspace{0.2cm} 0 \leq t \leq \tau_c(\nu_x). &
\end{flalign*}
In particular, if the mean curvature $\has$ of $\gamze$ is constant along $\gamze$ then the solution $\textsl{y}$ described in Lemma \ref{lemma_comparing_ywithJandj} and the function $\mathbb{X}= \mathbb{X}(x,t) $ do not depend on $x$ (that will be the case, for instance, for slices in a class of warped product spaces that we will analyse in a moment). However, for the comparison (\ref{willmorex_type_inequality_main}) to be accomplished with $\gamze$ in place of $\Sigma$ we must further assume that $(\textsl{M}, \textsl{g})$ has a well defined asymptotic volume ratio. Based on this, we assume that:

\begin{enumerate}[label=$(\mathcal{VR})$,leftmargin=1.7cm]
\label{asp_vol_rti_str}
\item the function $(0, \infty) \ni R \mapsto \Theta^{\text{\ding{93}}}(R) \doteq \frac{\textit{vol}(\mathcal{T}_{\gamze}^{+}(R))}{\w _{n} R^{n}}$ is (eventually) nonincreasing.
\end{enumerate}
Then, define the asymptotic$^{\text{\ding{93}}}$ volume ratio of $\textsl{g}$ as

\begin{flalign*}
\hspace{0.8cm}
AVR^{\text{\ding{93}}}(\textsl{g}) & \doteq \lim_{R \rightarrow \infty} \Theta^{\text{\ding{93}}}(R). &
\end{flalign*}
We should think of $\textit{vol}(\mathcal{T}_{\gamze}^{+}(R))$ as the volume of an oriented tubular neighborhood, where the orientation was chosen relatively to some unbounded component, namely $\textsl{M}^+$. Observe that if $\gamze$ is null-homologous, so that $\1\Omega = \gamze$ for some open and bounded set $\Omega$ then $\textsl{M} \setminus \gamze = (\textsl{M} \setminus \bar{\Omega}) \cup \Omega$, with $\textsl{M} \setminus \bar{\Omega}$ occupying the role of the unbounded component $\textsl{M}^{+}$ and 

\begin{flalign*}
\hspace{0.8cm}
\lim_{R \rightarrow \infty} \frac{\textit{vol}(\mathcal{T}_{\gamze}^{+}(R))}{\w_n R^n} & = \lim_{R \rightarrow \infty} \frac{ \textit{vol}(\mathcal{T}_{\Omega}(R))}{\w_n R^n}&
\end{flalign*}
In consequence, the definition of $AVR^{\text{\ding{93}}}(\textsl{g})$ is natural, although not through geodesic balls. Proceeding exactly as in \hyperref[prf_ineq_four]{Proof of inequality (4)}, we easily obtain the inequality (\ref{wil_typ_inq_sec_main}) for $\gamze$ instead of $\Sigma$ under these conditions. In particular, this exposition implies our inequality (\ref{willmorex_type_inequality_wp}).$\hspace{0.3cm} \square$ \vspace{0.2cm}
\end{dig}

Next, we show how the proof of Theorem \ref{main_theo} may be adapted to the class of asymptotic$^{\text{\ding{93}}}$ nonnegatively curved spaces under the circumstances described above. \vspace{0.2cm}

\textit{Proof of Theorem \ref{sec_main_theo}.} Here we use both the distances $w \mapsto\text{dist}_{\gamze}(w)$ and $r^{\bullet}(w) = \text{dist}_{\Sigma}(w)$, where $w \in \textsl{M}^{+}$, and the proof follows virtually the same lines of Theorem $\ref{main_theo}$ by switching the distance $\distP$ by $\text{dist}_{\gamze}$, so we skip minor details. We must analyse the asymptotics of

\begin{flalign*}
\hspace{0.8cm}
\frac{\textit{vol}(\mathcal{T}_{\Sigma}^{+}(R))}{\vert \mathbb{B}^{n}(R) \vert} & = \frac{n}{R^n \vert \mathbb{S}^{n-1} \vert}  \int_{\Sigma} \int_{0}^{R \wedge \tau_c(x)} \dete J_{\xi(x)} (t) dt d\sigma (x), &
\end{flalign*}
where $\tau_c(x)$ refers to the cut time of the geodesics with initial velocity $\xi_x$. For $x \in \Sigma$, define the function $\mathscr{J}_{x}^{\text{\ding{93}}} : [0, t_{\textsl{f}}(\xi_x)) \to \mathbb{R}$ by

\begin{flalign}
\label{def_J_str}
\hspace{0.8cm}
\mathscr{J}_{x}^{\text{\ding{93}}}(t) &\doteq \big(\dete J_{\xi(x)}(t)\big)^{1/(n-1)}.&
\end{flalign}
Since $\Ricc_w \geq  -(n-1) (\lambda \circ \text{dist}_{\gamze})(w)\textsl{g}_w$, and $\lambda \neq 0$ proceeding as in Section \ref{section_elel_ineq} we readily obtain

\begin{itemize}
\item $(\mathscr{J}_{x}^{\text{\ding{93}}})^{\prime \prime}(t) - \big(\lambda \circ \text{dist}_{\gamze} \circ \gamma_{\xi_x}\big)(t) \mathscr{J}_{x}^{\text{\ding{93}}}(t)  \leq 0,$ and
\item $\mathscr{J}_{x}^{\text{\ding{93}}} (t)\leq \textsl{y}_{x}^{\text{\ding{93}}}(t) \leq \mathbb{X}_{x}(t),$ 
\end{itemize}
for all $x \in \Sigma$ and all $t \in [0, t_{\textsl{f}}(\xi_x))$, where the function $\mathbb{X}_{x}$ is defined in (\ref{definition_of_X_f}) and $\textsl{y}_{x}^{\text{\ding{93}}}$ is the unique solution to the following $IVP$

\begin{alignat*}{2}
\hspace{-4.5cm}
   & \begin{aligned} & \begin{cases}
  (\textsl{y}_{x}^{\text{\ding{93}}})^{\prime \prime}(t) - \lambda (\text{dist}_{\gamze} (\gamma_{\xi_x}(t)))\textsl{y}_{x}^{\text{\ding{93}}}(t) =0, & \\
 	\textsl{y}_{x}^{\text{\ding{93}}}(0)=1, \hspace{0.2cm} (\textsl{y}_{x}^{\text{\ding{93}}})^{\prime}(0)=\vert \hit \vert(x), \hspace{1.5cm} & \\
  \end{cases}\\
%  \MoveEqLeft[-1]\text{}
  \end{aligned}
\end{alignat*}
where $\hit = \has/(n-1)$, with $\has$ being the mean curvature of $\Sigma$. By hypothesis, there exists an open set $\Omega$ with finite volume such that $\1 \Omega = \gamze \cup \Sigma$. Then $\textit{vol}(\mathcal{T}_{\gamze}^{+}(R)) \leq \vert \Omega \vert + \textit{vol}(\mathcal{T}_{\Sigma}^{+}(R))$ so that, proceeding exactly as in \hyperref[prf_ineq_four]{Proof of inequality (4)} and letting $R \rightarrow \infty$ in a chain of inequalities of the type

\begin{flalign*}
\hspace{0.8cm}
\frac{\textit{vol}(\mathcal{T}_{\gamze}^{+}(R))}{\w_n R^n}
	& \leq \frac{\vert\Omega\vert}{\vert \mathbb{B}^{n}(R) \vert} + \frac{1}{\vert \mathbb{B}^{n}(R)\vert}  \int_{\Sigma} \int_{0}^{R \wedge \tau_c(x)} \mathscr{J}^{\text{\ding{93}}}(x,t)^{n-1} dt d\sigma (x) \\
	& \leq  \frac{\vert\Omega\vert}{\vert \mathbb{B}^{n}(R) \vert} + \frac{n}{R^n \vert \mathbb{S}^{n-1} \vert}  \int_{\Sigma} \int_{0}^{R} \mathbb{X}(x,t)^{n-1} dt d\sigma (x) &
\end{flalign*}
grants us inequality (\ref{wil_typ_inq_sec_main}). \vspace{0.2cm}

\textbf{Equality discussion.} Suppose we have the equality

\begin{flalign}
\label{eq_asym_star}
\hspace{0.8cm}
e^{(n-1)b_0}\int_{\Sigma}\bigg(\bigg \vert \frac{ \has}{n-1} \bigg \vert (1+b_0) + b_1 \bigg)^{n-1}d\sigma  & = AVR^{\text{\ding{93}}}(\textsl{g}) \vert\mathbb{S}^{n-1} \vert. &
\end{flalign}
First observe $\tau_c \equiv \infty$ over $\Sigma$ along the directions given by $\xi$. Analogously to what was done in Digression \ref{dig_on_collective_behavior} we define the functions $\hat{\theta}^{\text{\ding{93}}}, \theta^{\text{\ding{93}}}: \Sigma \times [0, \infty) \to \mathbb{R}$ setting

\begin{flalign*}
\hspace{0.8cm} \hat{\theta}^{\text{\ding{93}}}(x, t)^{\frac{1}{n-1}} & = \bigg(\frac{\mathscr{J}^{\text{\ding{93}}}(x, t)}{\mathbb{X}(x,t)} \bigg) \hspace{0.2cm} \text{and}\hspace{0.2cm}
\theta^{\text{\ding{93}}}(x, t)^{\frac{1}{n-1}}= \bigg(\frac{\textsl{y}^{\text{\ding{93}}}(x, t)}{\mathbb{X}(x,t)}\bigg).&
\end{flalign*}
Now, the fact that $\frac{\mathscr{J}^{\text{\ding{93}}}(x,t)}{\mathbb{X}(x,t)} \leq \frac{\mathscr{J}^{\text{\ding{93}}}(x,t)}{\textsl{y}^{\text{\ding{93}}}(x,t)}$, $x \in \Sigma$, $t \geq 0$ together with the monotonicity of $t \mapsto \frac{\mathscr{J}_x^{\text{\ding{93}}}(t)}{\textsl{y}_{x}^{\text{\ding{93}}}(t)}$ and equality (\ref{eq_asym_star}) gives $\sup\{ \hat{\theta}^{\text{\ding{93}}}(x,t): t \in [R^{\prime}, \infty) \} = 1$, for almost every $x \in \Sigma$ and for all $R^{\prime}>0$ which yields $\mathscr{J}^{\text{\ding{93}}}(x,t) \equiv \textsl{y}^{\text{\ding{93}}}(x,t)$ on $\Sigma \times [0, \infty)$ so that the mean curvature $\has$ of $\Sigma$ is nonnegative and 

\begin{enumerate}[label=\textit{(E\hspace{-0.035cm}\arabic*$^{\text{\ding{93}}}$)},leftmargin=1.7cm]
\item $\mathcal{H}_{r^{\bullet}} = \frac{\Delta r^{\bullet}}{n-1} \mathcal{P}_{\hspace{-0.08cm}\1_{r^{\bullet}}}$;

\item $\Ricc \Big( \frac{\1}{\1 r^{\bullet}} , \frac{\1}{\1 r^{\bullet}}\Big)=-(n-1)(\lambda \circ \text{dist}_{\gamze})$,
\end{enumerate}
both hold on $\Phi(\Sigma \times[0, \infty) )$, where $\Phi$ is defined by (\ref{diffeo_Phi}). Repeating the argument in the \textbf{Equality discussion} in Section \ref{sharp_dom_asym_non} we deduce that $\Sigma$ is a totally umbilic hypersurface with locally constant mean curvature. In addition to that, if $\Sigma$ is connected then \textit{(E\hspace{-0.035cm}1$^{\text{\ding{93}}}$)} and \textit{(E\hspace{-0.035cm}2$^{\text{\ding{93}}}$)} together with the nonexistence of critical points for $r^{\bullet}$, the fact $\tau_c \equiv \infty$ over $\Sigma$ and elementary Morse theory yield that the referred distance induces a smooth, codimension one umbilical foliation $\mathcal{F}^{\text{\ding{93}}}(\textsl{X})$ whose leaves $\Sigma_{t}^{\text{\ding{93}}} = \Phi \big(\Sigma \times \{ t \}  \big)$ are connected hypersurfaces with constant mean curvature $\has[t]$, whence $\textsl{M}^{+} \setminus \Omega$ has only one end. The aforementioned foliation is also induced by the closed and conformal vector field given by

\begin{flalign*}
\hspace{0.8cm}
\textsl{X}\big\vert_{w} & = \varphi(r^{\bullet}(w)) \1_{r^{\bullet}}\big\vert_{w}, \hspace{0.2cm}\text{and defined on} \hspace{0.2cm} \textsl{M}^{+} \setminus \Omega \hspace{0.2cm} \big( \text{observe} \hspace{0.2cm} C(\Sigma) = \emptyset \big) &
\end{flalign*}
where $\varphi$ is the solution to $\varphi^{\prime}(t)- \frac{\has[t]}{n-1}\varphi(t)=0$, with initial condition $\varphi(0)=1$. Following the ideas of Montiel (see Proposition 3 in \cite{montiel1998stable}) one shows that $\textsl{M}^{+} \setminus \Omega$ is locally a warped product. Following the ideas in our \textbf{Equality discussion} in the preceding section, one shows that

\begin{flalign*}
\hspace{0.8cm} &
\textsl{M}^{+} \setminus \Omega \hspace{0.2cm} \text{is isometric to} \hspace{0.2cm} \big([r_0, \infty)\times \Sigma, dr \otimes dr + \varrho(r)^2 \textsl{g}_{\Sigma}\big), \hspace{0.2cm}\text{where} \hspace{0.2cm} r_0 = \text{dist}_{\textsl{g}}(\Gamma_0,\Sigma)&
\end{flalign*}
and $\varrho(r)= \varphi(r-r_0)$ is the solution to $\varrho^{\prime \prime}(t) - \lambda(t)\varrho(t)=0$, with the required initial conditions, so that $\gamze$ and $\Sigma$ are equidistant hypersurfaces (see Remark \ref{obs_sup_varphi}) with $r_0$ being the distance between them. $\hspace{0.3cm} \square$ \vspace{0.2cm}

\section{Application to warped product spaces}
\label{app_und_exm}

In this section, we show that Corollary \ref{main_prop} is a direct consequence of Theorem \ref{sec_main_theo} and set the stage for the proof of the converse of our Willmore-type inequality in asymptotic and asymptotic$^{\text{\ding{93}}}$ nonnegatively curved spaces. Moreover, we evaluate explicitly the constants $b_0(\lambda)$, $b_1(\lambda)$ and $AVR^{\text{\ding{93}}}(\textsl{g}_{\textit{S}})$ for the 3-dimensional Schwarzschild manifold. \vspace{0.2cm}

\begin{lemma} $\big[$\textit{Basic properties of $(\textsl{M}, \textsl{g})$} $\big]$ \newline
\label{bas_prp_mg} 
Let $(\textsl{M}, \textsl{g})$ be the warped product described in the introduction, that is

\begin{flalign*}
\hspace{0.8cm}
\textsl{M} & \doteq [0, \infty)\times \textsl{N}, \hspace{0.2cm}
\textsl{g} \doteq dr \otimes dr + h(r)^2 \textsl{g}_{N}, &
\end{flalign*}
and satisfying the conditions $(\mathbf{\Lambda}\hspace{-0.02cm}\textbf{1})\textit{---}(\mathbf{\Lambda}\hspace{-0.02cm}\textbf{3})$. Then, $(\textsl{M}, \textsl{g})$ has asymptotically$^{\text{\ding{93}}}$ nonnegative Ricci curvature relatively to $\gamze$ and associated function $\lambda$ with finite $AVR^{\text{\ding{93}}}(\textsl{g})$, as defined in Digression \ref{dig_hom_hyp}. In addition to that, if the condition $(\mathbf{\Lambda}\hspace{-0.02cm}\textbf{4})$ holds then $AVR^{\text{\ding{93}}}(\textsl{g})>0.$
\end{lemma}
\textit{Proof.}
Let $\{ E_{1} \doteq \1_r,E_2 ..., E_{n} \}$ be a local orthonormal frame for $\textsl{M}$ so that $\textsl{g}(E_l, E_{k}) = \delta_{lk}$ and set $e_{k} = hE_k$, for $k=2,...,n$. It follows from elementary properties of warped product spaces (see e.g. \cite{gromoll2009metric}, pp. 59 \textit{--} 60) that the Ricci curvature of $\textsl{M}$ in the radial direction is given by 

\begin{flalign*}
\hspace{0.8cm}
\Ricc_{\textsl{M}}(\1_r, \1_r) & = -(n-1)\frac{h^{\prime \prime}(r)}{h(r)}. &
\end{flalign*}
Whereas in spacial directions

\begin{flalign*}
\hspace{0.8cm}
\Ricc_{\textsl{M}}(E_k, E_k) & = \Ricc_{\textsl{N}}(E_k, E_k)-
\bigg( \frac{h^{\prime \prime}(r)}{h(r)} + (n-2)\frac{h^{\prime}(r)^2}{h(r)^2}\bigg) \notag \\
	& \geq -(n-1) \bigg[ \frac{1}{(n-1)} \frac{h^{\prime \prime}(r)}{h(r)} - 
	\bigg(\frac{n-2}{n-1}\bigg) \frac{\rho - h^{\prime}(r)^2}{h(r)^2}\bigg] . &
\end{flalign*}
Let $\lambda$ be the function from hypothesis $(\mathbf{\Lambda}\hspace{-0.02cm}\textbf{2})$ so as to $\lambda$ is nonincreasing and $b_0$ is well defined. The computations above show that

\begin{flalign*}
\hspace{0.8cm}
\Ricc_{\textsl{M}}& \geq  -(n-1)\big( \lambda \circ \text{dist}_{\gamze}\big)\textsl{g} &
\end{flalign*}
so that the first claim is verified. If $(\mathbf{\Lambda}\hspace{-0.02cm}\textbf{3})$ holds then $t \mapsto \frac{h(t)^{n-1}}{t^{n-1}}$ is eventually nonincreasing. By Lemma 2.2 in \cite{10.2307/2374996} the function

\begin{flalign*}
\hspace{0.8cm}
(0, \infty) \ni R \hspace{0.1cm} \mapsto & \hspace{0.15cm} \frac{1}{{\displaystyle \int_{0}^{R}t^{n-1}dt }} \int_{0}^{R} h(t)^{n-1}dt &
\end{flalign*}
is also eventually nonincreasing. Since the (open) geodesic tube around $\gamze$ of radius $R < \infty$ is equal to $\mathcal{T}_{\gamze}(R) = [0, R) \times \textsl{N}$ we have

\begin{flalign}
\label{smp_vol_str_tub_exp}
\hspace{0.8cm}
\frac{\textit{vol}(\mathcal{T}_{\gamze}(R))}{\vert \mathbb{B}^{n}(R) \vert} & = 
\frac{\vert \textsl{N} \vert}{\vert \mathbb{S}^{n-1} \vert} \frac{1}{{\displaystyle \int_{0}^{R}t^{n-1}dt }} \int_{0}^{R} h(t)^{n-1}dt. &
\end{flalign}
Hence, for big values of radius $R$, the function $R \mapsto  \textit{vol}(\mathcal{T}_{\gamze}(R))\big/\vert \mathbb{B}^{n}(R)\vert$ is nonincreasing and we have a finite asymptotic$^{\text{\ding{93}}}$ volume ratio. If $(\mathbf{\Lambda}\hspace{-0.02cm}\textbf{4})$ holds then it is easy to see that $AVR^{\text{\ding{93}}}(\textsl{g}) >0$. $\hspace{0.3cm} \square$ \vspace{0.2cm}

We now reveal how Corollary \ref{main_prop} is obtained from Theorem \ref{sec_main_theo}. Choose the open set $\textsl{M}^{+}$ from Theorem \ref{sec_main_theo} to be equal to $\textsl{M} \setminus \gamze$ and $\textsl{M}^{-} = \emptyset$. Since $\gamze$ and $\Sigma$ are homologous, there exists an open set $\Omega$ with $\vert \Omega \vert <\infty$ such that $\1 \Omega = \Sigma \cup \gamze$. The second inequality in the corollary follows from here and the first was observed in Digression \ref{dig_hom_hyp}. With respect to the rigidity statement, if equality in (\ref{will_type_ine_connected}) holds then by Theorem \ref{sec_main_theo} the manifold $\textsl{M}\setminus (\Omega \cup \gamze)$ is isometric to the model (\ref{iso_asy_star_to}) (here we use the connectedness assumption on $\Sigma$), the hypersurfaces $\gamze$ and $\Sigma$ are equidistant and $\therefore$ $\Sigma$ is a slice $\{r_0\}\times \textsl{N}$, where $r_0$ is the distance between $\gamze$ and $\Sigma$. Clearly, $h \equiv \varrho$ over $[r_0, \infty)$ so that all the necessary conditions are met. The proof of the sufficiency along with the converse of Theorems \ref{main_theo} and \ref{sec_main_theo} will be shown in the Subsection \ref{converse}. \vspace{0.2cm}

At this point, considering the warped product from Lemma \ref{bas_prp_mg} it is interesting to look at the behavior of the function $\mathscr{F}:[0, \infty) \to \mathbb{R}$ given by

\begin{flalign*}
\hspace{0.8cm}
\mathscr{F}(t)& = \int_{\textsl{N}} \bigg(\bigg \vert \frac{ h^{\prime}(t)}{h(t)}\bigg \vert \big(1+b_0\big) + b_1 \bigg)^{n-1}\textit{d} \textit{vol}_{\textsl{N}}.&
\end{flalign*}
The slices $\{r\} \times \textsl{N}$ for which equality in (\ref{will_type_ine_connected}) occurs represent (global) minima of $\mathscr{F}$. If $(\mathbf{\Lambda}\hspace{-0.02cm}\textbf{4})$ holds so that $h^{\prime}(t) \geq 0$ for $t \geq \tau_0$ then

\begin{flalign*}
\hspace{0.8cm}
\mathscr{F}^{\prime}(t)& = \frac{(n-1)\vert \textsl{N} \vert}{(1+b_0)^{-1}} \bigg( \frac{ h^{\prime}(t)}{h(t)} \big(1+b_0\big) + b_1 \bigg)^{n-2} \bigg( \frac{ h^{\prime \prime}(t)}{h(t)} - \frac{ h^{\prime}(t)^2}{h(t)^2}\bigg), \hspace{0.2cm} t \geq \tau_0 .&
\end{flalign*}
If equality in (\ref{will_type_ine_connected}) holds for the slice $\Sigma = \{\tau_0\}\times \textsl{N}$ then $\mathscr{F}^{\prime}(t) \geq 0$ on a neighborhood of $\tau_0$ and

\begin{flalign*}
\hspace{0.8cm}
\lim_{t \rightarrow \infty} \mathscr{F}^{\prime}(t) & = 0. & 
\end{flalign*}
\vspace{0.2cm}

\subsection{The converse of Theorems \ref{main_theo} and \ref{sec_main_theo}}
\label{converse}

We are now ready to prove the converse of Theorems \ref{main_theo} and \ref{sec_main_theo}. If $\Omega$ is null-homologous (resp. homologous to $\gamze$) then $\textsl{M} \setminus \Omega$ (resp. $\textsl{M}^{+}\setminus \Omega$) is isometric to

\begin{flalign*}
\hspace{0.8cm} \biggl( [r_0, \infty) \times \Sigma & , dr\otimes dr + \varrho(r)^2 \textsl{g}_{\Sigma}\biggl) &
\end{flalign*}
where $r_0$ is the distance between $\pzero$ and $\Sigma$ (resp. $\gamze$ and $\Sigma$), $\Sigma$ is a connected constant mean curvature hypersurface and $\varrho$ is the solution to the Jacobi equation $\varrho^{\prime \prime} - \lambda\varrho =0$ with the given initial conditions and satisfying (\ref{iso_asy_to}) (resp. (\ref{iso_asy_to_star})). Now, for positive $R$ (resp. $R > r_0$) we have 
$\textit{vol}(\mathcal{T}_{\Omega}(R)) = \vert \Omega \vert + \textit{vol}(\mathcal{T}_{\Sigma}^{+}(R))$ (resp. $\textit{vol}(\mathcal{T}_{\gamze}^{+}(R)) = \vert \Omega \vert + \textit{vol}(\mathcal{T}_{\Sigma}^{+}(R-r_0)))$. Due to legibility, we now restrict ourselves to calculating $AVR(\textsl{g})$ and note that similar calculations hold for $AVR^{\text{\ding{93}}}(\textsl{g})$. According to formula (\ref{smp_vol_str_tub_exp}), the definition of the asymptotic volume ratio and the fact that $\mathcal{T}_{\Sigma}^{+}(R) = [r_0, r_0+R)\times \Sigma$, we have

\begin{flalign*}
\hspace{0.8cm}
AVR(\textsl{g})
	& = \lim_{R \rightarrow \infty} \Bigg( \frac{\vert \Sigma \vert}{\vert\mathbb{S}^{n-1}\vert} \bigg(\int_{0}^{R}t^{n-1}dt \bigg)^{-1} \int_{r_0}^{r_0+R} \varrho(t)^{n-1}dt \Bigg) &
\end{flalign*}
thereby

\begin{flalign*}
\hspace{0.8cm} 	
\frac{AVR(\textsl{g})\vert\mathbb{S}^{n-1}\vert}{\vert \Sigma \vert} 
	& =  \lim_{R \rightarrow \infty} 
\bigg( \frac{\varrho(r_0+R)}{R} \bigg)^{n-1} \\
	& = \bigg( \lim_{R \rightarrow \infty} \frac{\varrho(r_0+ R)}{R} \bigg)^{n-1} \\
	& = \bigg[ e^{b_0}\bigg(\frac{ \has }{(n-1)}(1+b_0) + b_1 \bigg)\bigg]^{n-1}\hspace{-0.09cm}. &
\end{flalign*}
This implies that we have equality in (\ref{willmorex_type_inequality_main}).

\subsection{Application to the Schwarzschild and Reissner-Nordstrom manifolds}

A particular set of examples of the above construction comes from considering Riemannian warped products $(\textsl{M}, \textsl{g})$ of the following type 

\begin{flalign}
\label{non_wp_metric}
\hspace{0.8cm}
 \textsl{M} & = (\sht, \infty) \times \textsl{N} , \hspace{0.2cm} \textsl{g} = \frac{1}{\w(s)}ds \otimes ds + s^2 \textsl{g}_{\textsl{N}},&
\end{flalign}
where $\sht >0$ and $\w(s)$ is a smooth function on $[\sht, \infty)$. As described in \cite{brendle2013constant}, to see the metric $\textsl{g}$ into the form (\ref{warped_product_metric}), one considers the change of variables $F: [ \sht , \infty) \to \mathbb{R}$ defined by $F^{\prime}(s)=1/\sqrt{\w (s)}$ and $F(\sht)=0$. Then, the substitution $r = F(s)$ brings (\ref{non_wp_metric}) into the desired form and if $h:[0, \infty)\to [\sht, \infty)$ denotes the inverse of $F$ an easy computation gives 

\begin{flalign*}
\hspace{0.8cm}
h^{\prime}(r)& = \sqrt{\w(s)} \hspace{0.2cm} \text{and} \hspace{0.2cm} h^{\prime \prime}(r)= \frac{1}{2}\w^{\prime}(s), &
\end{flalign*}
where $s=h(r)$. With respect to the change of variables $r=F(s)$, the functions $\lambda_1$ and $\lambda_2$ become

\begin{flalign*}
\hspace{0.8cm}
\lambda_1(s) & = \frac{\w^{\prime}(s)}{2s} \hspace{0.2cm} \text{and} \hspace{0.2cm} \lambda_2(s) = \frac{\w^{\prime}(s)}{2(n-1)s} - \bigg(\frac{n-2}{n-1}\bigg)\frac{\rho - \w(s)}{s^2}, &
\end{flalign*}
the condition ($\mathbf{\Lambda}$\hspace{-0.02cm}\textbf{1}) remains unchanged while ($\mathbf{\Lambda}$\hspace{-0.02cm}\textbf{2}) \textit{---} ($\mathbf{\Lambda}$\hspace{-0.02cm}\textbf{4}) are equivalent to the following pair of conditions

\begin{itemize}
\item there exists a continuous nonnegative nonincreasing function $\lambda: [\sht, \infty) \to \mathbb{R}$ s.t. $\lambda \geq \text{max}\{\lambda_1, \lambda_2\}$ and 

\begin{flalign}
\hspace{0.8cm}
b_0 & = {\displaystyle \int_{\sht}^{\infty}F(s)\lambda(s)F^{\prime}(s)ds} < \infty  &
\end{flalign}

\item the function $(\sht, \infty) \ni s \mapsto \frac{s}{F(s)}$ is eventually monotone decreasing.
\end{itemize}

It is straightforward to check that, after the change described above, the Schwarzschild and Reissner-Nordstrom spaces introduced in (\ref{n-d_sch_mtr}) and (\ref{n-d_rei_nord_mtr}) respectively, fulfils the conditions $(\mathbf{\Lambda}\hspace{-0.02cm}\textbf{1})\textit{---}(\mathbf{\Lambda}\hspace{-0.02cm}\textbf{4})$. In particular, since these spaces have smallest Ricci curvature in the radial direction, there is an evident choice for the function $\lambda$. \vspace{0.2cm}

\begin{exmp}
\label{ex_schwarzs} $\big[$\textit{The constants $b_0(\lambda)$, $b_1(\lambda)$ and $AVR^{\text{\ding{93}}}(\textsl{g}_{\textit{S}})$ for the 3-d Schwarzschild space} $\big]$ \newline
The (3-dimensional) model of the Schwarzschild space we are adopting in this example is 

\begin{flalign}
\label{3-d_sch_mtr}
\hspace{0.8cm}
 \textsl{M} & = [0, \infty) \times \mathbb{S}^{2}, \hspace{0.2cm} \textsl{g} = dr \otimes dr + h(r)^2 \textsl{g}_{\mathbb{S}^{2}},&
\end{flalign}
where $h$ has $h^{\prime}(r) = \sqrt{1-ms(r)^{-1}}$ so that it meets the initial condition $h(0)=m$. It is straightforward to check that $\textsl{M}$ has asymptotically$^{\text{\ding{93}}}$ nonnegative Ricci curvature relatively to $\gamze = \{0\} \times \mathbb{S}^2$:
\begin{flalign*}
\hspace{0.8cm}
\Ricc_{\textsl{M}} & \geq -(n-1)\frac{h^{\prime \prime}(r)}{h(r)}\textsl{g}. &
\end{flalign*}
In this case we have the clear choice $\lambda (r) = \frac{h^{\prime \prime}(r)}{h(r)}$ for the associated function. Also, in this example, $\Sigma$ is taken to be equal to $\gamze$. We begin by computing the constant $b_1$. Using the substitution $r=F(s)$, the definition of $b_1$ and the change of variables $u = \w(s)$, where $\w(s) = 1-ms^{-1}$, we get

\begin{flalign*}
\hspace{0.8cm}
b_1 & = \int_{0}^{\infty} \lambda(t)dt = \frac{m}{2} \int_{\sht}^{\infty} \frac{s^{-3}}{\sqrt{\w(s)}} ds=\frac{2}{3} \frac{1}{m} &
\end{flalign*}
We now evaluate $AVR^{\text{\ding{93}}}(\textsl{g})$. Using formula (\ref{smp_vol_str_tub_exp}), using the substitution $r=F(s)$ and the change of variables $u = \w(s)$ again, we have

\begin{flalign*}
\hspace{0.8cm}
\textit{vol}(\mathcal{T}_{\Sigma}(R)) & = 
\frac{m^3 \vert \mathbb{S}^2 \vert}{48} \Bigg\{ \frac{2\sqrt{\w}}{m^{3}s(R)^{-3}}\big( 15\w ^2 -40\w +33\big) -15 \lna \bigg( \frac{ms(R)^{-1}}{(1+\sqrt{w})^2} \bigg) \Bigg\}, \hspace{0.2cm} \w = \w(s(R)).&
\end{flalign*}
Therefore,

\begin{flalign*}
\hspace{0.8cm}
AVR^{\text{\ding{93}}}(\textsl{g})
	& = \lim_{R\to \infty}\frac{\textit{vol}(\mathcal{T}_{\Sigma}(R))}{\vert \mathbb{B}^{3}(R) \vert} \\
	& = \lim_{R\to \infty} \frac{m^3}{8 R^3} \Bigg\{ \frac{\sqrt{\w}}{m^{3}s(R)^{-3}}\big(15\w^2 -40\w +33\big) -\frac{15}{2} \lna \bigg( \frac{ms(R)^{-1}}{(1+\sqrt{w})^2} \bigg) \Bigg\} =1.&
\end{flalign*}
To compute $b_0(\lambda)$ note that $F(s) =(1/2) \big(2s\sqrt{\w}+ m\lna(1+ \sqrt{\w}) - m\lna (1- \sqrt{\w})\big)$.
Then,

\begin{flalign}
\label{b_zero}
\hspace{0.8cm}
b_0 & = \frac{m}{2}\int_{\sht}^{\infty} \frac{F(s)}{s^3\sqrt{\w(s)}}ds = \frac{1}{3}\big(1+ \lna 4 \big). \hspace{0.3cm} \square \vspace{0.2cm}& 
\end{flalign}
\end{exmp}

\vspace{0.2cm}

Observe that it is natural to make an attempt to extract a model for the geometry (\ref{iso_asy_star_to}) from the Schwarzschild space by changing the warping function. We now describe one such experiment. \vspace{0.2cm}

Let $\kappa$ be a positive real number and define $m \doteq \frac{2}{3\kappa} e^{b_0}$, where $b_0$ is given by (\ref{b_zero}). As in the previous example, choose the warping function $h$ with initial condition $h(0)= m$. Then define $\varrho:[0, \infty) \to \mathbb{R}$ by $\varrho(r) \doteq m^{-1} h(r)$, so that $\varrho(0)=1$ and $\varrho^{\prime}(0)=0$
and consider the Riemannian manifold 

\begin{flalign*}
%\label{3-d_sch_mtr}
\hspace{0.8cm}
 \textsl{M} & = [0, \infty) \times \mathbb{S}^{2}, \hspace{0.2cm} \textsl{g} = dr \otimes dr + \varrho(r)^2 \textsl{g}_{\mathbb{S}^{2}}.&
\end{flalign*}
It is straightforward to check that the Ricci curvature of $\textsl{g}$ is smallest in the radial direction provided $m \geq 1$. Therefore, $\Ricc_{\textsl{M}} \geq -(n-1)(\lambda \circ \text{dist}_{\gamze})\textsl{g}$, where the associated function $\lambda$ is equal to $\varrho^{\prime \prime}/\varrho = h^{\prime \prime}/h$ so that it satisfies the Jacobi equation $\varrho^{\prime \prime} - \lambda \varrho = 0$. As a result, the numbers $b_1$ and $b_0$ remain the same as those calculated in the previous example. However,

\begin{flalign*}
\hspace{0.8cm} & \lim_{r \rightarrow \infty}\frac{\varrho(r)}{\big(e^{b_{0}} b_1\big)r} = \lim_{r \rightarrow \infty} \frac{3}{2e^{b_0}} \frac{h(r)}{r}  = \frac{3}{2e^{b_0}} \neq 1, &
\end{flalign*}
and we do not achieve equality in (\ref{willmorex_type_inequality_wp}).
\vspace{0.2cm}

\section{Annexes}
\label{zeta}
In this section we discuss in more depth some topics left along the text. Firstly, regarding the Willmore-type inequality in spaces of nonnegative Ricci curvature, we provide sufficient conditions on the open set $\Omega$ that guarantees that the whole space is isometric with the Euclidean space. Secondly, we show how to obtain a weaker version of the Willmore-type inequality using Sobolev inequalities, as discussed in the introduction.

\subsection{The geometry of $\Omega$}
\label{geometry_of_Omega}

Here we consider the equality case in equation (\ref{willmore_type_inequality_main_ricci_nonnegative}), i.e., suppose

\begin{flalign*}
\hspace{0.8cm}
\int_{\Sigma}\bigg \vert \frac{ \has}{n-1} \bigg \vert ^{n-1}d\sigma  & = AVR(\textsl{g}) \vert\mathbb{S}^{n-1} \vert >0.&
\end{flalign*}
Then, as described in \cite{agostiniani2020sharp}, $(\textsl{M},\textsl{g})$ has Euclidean Volume Growth and $\textsl{M}\setminus \Omega$ may be written as the warped product $[r_0, \infty) \times_{h} \Sigma$ endowed with metric $dr\otimes dr+h(r)^{2}\textsl{g}_{\Sigma}$, where $h(r)=r/r_0$ and $\Sigma$ is a connected totally umbilic hypersurface with constant mean curvature $\has$. The nonnegativity of the Ricci tensor together with elementary properties of warped product spaces (cf. Corollary 2.2.2. in \cite{gromoll2009metric}) gives

\begin{flalign*}
\hspace{0.8cm}
\Ricc_{\Sigma} & \geq  (n-2)\frac{1}{r_{0}^{2}} \textsl{g}_{\Sigma }. &
\end{flalign*}
It then becomes noteworthy to ask what geometric properties of $\Omega$, or its boundary, prevent the whole space to be isometric with $\mathbb{R}^{n}$ with standard Euclidean metric. It is straightforward to check that if the diameter or the area of $\Sigma$ fulfills, respectively, $\text{diam}_{\textsl{g}}(\Sigma) \geq  \pi r_0$ or $\vert \Sigma \vert \geq r_{0}^{n-1} \vert \mathbb{S}^{n-1} \vert$ then, in fact, $(\textsl{M},\textsl{g})$ is isometric to the Euclidean space. An open question is if $\textit{vol}(\Omega) = \w_n$ yields the same conclusion. \hspace{0.2cm}

In a different direction since $\bar{\Omega}$ is compact with MCB we may use a sharp result for $\bar{\Omega}$ based only on its intrinsic distance to yield global conclusions. More specifically, by a theorem of Man-chun Li \cite{li2014sharp}, since $\Sigma$ has constant mean curvature, say $\has =(n-1)\textsl{k} >0$, we have 

\begin{flalign}
\label{equal_m_chun_li}
\hspace{0.8cm}
\sup_{x \in \Omega} \text{dist}_{\textsl{g}}(x, \Sigma) & \leq \frac{1}{\textsl{k}}. &
\end{flalign}
Furthermore, equality in (\ref{equal_m_chun_li}) holds if and only if $\bar{\Omega}$ is isometric to a closed Euclidean ball of radius $1/\textsl{k}$, in which case $AVR(\textsl{g}) = 1$ so that $(\textsl{M}, \textsl{g})$ is isometric to the Euclidean space.

\subsection{A weaker Willmore-type inequality}
\label{weaker_willmore_type}

We now describe how to obtain a weaker version of the Willmore-type inequality using Sobolev inequalities in spaces of nonnegative sectional curvature. Let $\Sigma^{n-1} \subset \textsl{M}$ be any closed hypersurface in a complete noncompact Riemannian manifold $(\textsl{M},\textsl{g})$ with nonnegative sectional curvature. Assume $\Sigma = \partial \Omega$, where $\Omega$ is open and bounded in $\textsl{M}$. Then, take $f=1$ over $\Sigma$ and, via the identification $\textsl{M} \cong \textsl{M} \times \{0\}$, realize $\Sigma$ as a codimension 2 submanifold in $\textsl{M}\times \mathbb{R}$, which is endowed with the product metric $\tilde{\textsl{g}}=\textsl{g}+dr\otimes dr$, apply Corollary 1.5 from \cite{brendle2023sobolev} to get

\begin{flalign}
\label{brendle_1.5_coroll}
\hspace{0.8cm}
\int_{\Sigma} \vert \vec{H} \vert d \sigma & \geq (n-1)\vert \mathbb{B}^{n-1} \vert^{\frac{1}{n-1}} \vert AVR(\tilde{\textsl{g}})^{\frac{1}{n-1}} \vert \Sigma \vert^{\frac{n-2}{n-1}}, &
\end{flalign}
where $\mathbb{B}^{k}$ is the euclidean unit ball in $\mathbb{R}^{k}$. A simple application of Hölder's inequality with $p=n-1$ and $q=(n-1)/(n-2)$ on the left hand side of the above yields

\begin{flalign*}
\hspace{0.8cm}
\bigg( \int_{\Sigma} \vert \vec{H} \vert^{n-1} d \sigma \bigg)^{\frac{1}{n-1}} & \geq (n-1) \vert \mathbb{B}^{n-1}\vert^{\frac{1}{n-1}} AVR(\textsl{g})^{\frac{1}{n-1}}, &
\end{flalign*}
where we have used $AVR(\textsl{g})= AVR(\tilde{\textsl{g}})$. As a result,

\begin{flalign*}
\hspace{0.8cm}
\int_{\Sigma} \bigg\vert \frac{\vec{H}}{n-1}  \bigg\vert^{n-1} d \sigma & \geq \frac{\vert \mathbb{S}^{n-2} \vert}{n-1} AVR(\textsl{g}). &
\end{flalign*}

\section{Acknowledgements}

This study was financed in part by the Coordenação de Aperfeiçoamento de Pessoal de Nível Superior – Brasil (CAPES) in terms of a doctoral scholarship. The author is very grateful to Professors Ezequiel Barbosa and Celso Viana for useful comments and discussions.

\bibliographystyle{unsrt}
\bibliography{references}
\end{document}